\documentclass[aps,pre,twocolumn,amsfonts,amsmath,showpacs,a4paper,floatfix]{revtex4-1}

\usepackage{enumerate}
\usepackage{amsfonts}
\usepackage{amssymb}
\usepackage{amsmath}
\usepackage{amscd}
\usepackage{psfrag}
\usepackage{latexsym}
\usepackage[dvips]{graphicx}
\usepackage{amsthm}
\usepackage{physics}
\usepackage{array}
\usepackage{float}
\usepackage{xcolor}

\newcommand{\seclabel}[1]{\label{sec:#1}}

\begin{document}

\title{Subcritical Hopf Bifurcations in the FitzHugh-Nagumo Model}

\author{S. Sehgal}
\email[]{sehgals@hope.ac.uk}
\affiliation{Department of Mathematics and Computer Science, Liverpool Hope University, Liverpool, UK}

\author{A. J. Foulkes}
\email[]{foulkea@hope.ac.uk}
\affiliation{Department of Mathematics and Computer Science, Liverpool Hope University, Liverpool, UK}

\date{\today}

\begin{abstract}
It had been shown that the transition from a rigidly rotating spiral wave to a meandering spiral wave
is via a Hopf bifurcation. Many studies have shown that these bifurcations are supercritical, but 
we present numerical studies which show that subcritical bifurcations are also present within FitzHugh-Nagumo.
Furthermore, we present a evidence that this bifurcation is highly sensitive to initial conditions, 
and it is possible to convert one solution in the hysteresis loop to the other.
\end{abstract}

\pacs{
	02.70.-c, 
  05.10.-a, 
  82.40.Bj, 
  82.40.Ck, 
  87.10.-e  
}

\maketitle

\section{Introduction}\seclabel{intro}

Spiral waves occur naturally in many physical, chemical and biological systems 
\cite{turing,wiener_ros, zhabot_1964a, zhabot_1964b, belousov_1959, FHN1, FHN2, 
Alcantara-Monk-1974, Oswald-Dequidt-2008, Dahlem-Mueller-2003}. 
The motion and behaviour of such waves can enlighten certain characteristics of the systems in 
which they occur. For instance, in cardiac tissue, the presence of these rotating spiral 
waves (also known as autowaves, and rotors)
indicates that there is an abnormality in the hearts natural rhythm 
(an arrhythmia) \cite{fenton98, keener_sneyd, moe77, winf72, pertsov93}. 
In most electro-chemical systems, such as cardiac tissue or neurological systems, 
excitable properties are an essential part of creating and sustaining spiral waves. 
The cells ability to be stimulated in response to external energy is critical
in the life cycle of spiral waves \cite{zykov18}. 

Excitable systems, such as the propagation of electrical energy along nerves, have been studied mathematically since 
1940 using parameter dependent mathematical models. Spiral waves were first observed by Wiener 
and Rossenblueth, who created the first finite automata model to simulate spiral 
wave activity \cite{wiener_ros}. Since then, many models of cell excitation
have been developed and studied. A classic system in modelling cellular excitation 
is the Hodgkin-Huxley model of nerve excitation developed in the early 
1950's \cite{hh52}. This model simulates the electrical energy passing through a 
single cell \cite{Schiening12}. It does not however simulate how each cell reacts with other cells 
that are part of the whole excitable medium. In order to explore this, spatial 
variables must be implemented and this is usually done using reaction-diffusion models.
Although the development of this type of model was initiated in the early part 
of the 20th century \cite{thomp17}, it was Turing in 1952 that used them to study 
interactions between chemical compounds \cite{turing}. Around the same time, 
Richard FitzHugh was developing a mathematical model to mimic threshold phenomena
in the nerve membrane \cite{fitz55}. In 1961, FitzHugh published a paper suggesting 
a model of nerve cell excitation, a simplified version of the Hodgkin-Huxley 
model, influenced by the Van der Pol oscillator equations \cite{FHN2}. 
In 1962, an equivalent model was publish by Nagumo et. al. \cite{FHN2}. 
The system of equations is now known as the FitzHugh-Nagumo model (FHN). To this 
day, many new models are developed using FHN as a base \cite{Bark91}. 

Mathematically, spiral waves are parameter dependent, spatio-temporal solutions 
to reaction-diffusion equations. Their motion relates directly to the physical 
phenomena that they represent. There are several main categories of motion of 
spiral wave, which are classified by considering how the \emph{tip}
of the spiral wave behaves and moves around the medium. 
The tip of the spiral wave can be defined as the intersection of two  
isolines in the excitation and inhibitor fields. If the tip traces out 
a perfect circle around a fixed center of rotation, then it is known as 
\emph{rigidly rotating}. A property of rigid rotation is that 
the shape of the arm of the spiral is fixed and the motion is periodic. 
Another type of motion is known as \emph{meander}. This is quasi-periodic with 
the arm of the spiral periodically changing shape and the tip of the spiral
tracing out epicycloid type patterns \cite{Winfree-1991,Wulff-1996,Barkley-1994,Barkley-1995}.

Other types of motion include hypermeander (chaotic motion) and drift 
(motion around a moving center of 
rotation) \cite{Ashwin-etal-2001,Biktashev-Holden-1998,Henry-2004,Biktashev-Holden-1995}. 
We are not considering these 
types of motion in this publication.

The transition between one type of motion to another via a parameter change 
has been of great interest since the study of spiral waves began. In particular, the
transition from rigid rotation to meander has been shown numerically to be via a 
Hopf bifurcation \cite{Bark90,bik99,Golub97}. In particular, several authors have noted that this 
transition is via a supercritical Hopf bifurcation, in which a stable fixed point
changes to unstable and a stable limit cycle is formed. A key feature of this
transition is that the growth of the limit cycle from the bifurcation point 
is proportional to the square root of the varying parameter.

Here we present work on a numerical approach to study the transition from 
rigid rotation (RW) to meander (MRW), and show that in the FHN system of equations, 
there are regions within the parameter space in which the Hopf bifurcations
are in fact subcritical. We analyse these results to confirm their validity 
and show that within the hysteresis regions where there are two solutions relating 
to the same set of parameters, it is possible to convert one of these solutions to the other.
\section{Numerical simulations of spiral waves in a comoving frame of reference}\seclabel{nums}
\subsection{Theoretical review}

Our approach to studying the transition from RW to MRW is motivated by the restrictions of 
earlier numerical studies in studying large core spirals \cite{Bark90}. The problems with large core 
spirals is that if we study them in a numerical box that is fixed in space, then that
box will need to be very large for spiral waves with large cores. This will result in very computationally
expensive simulations which will take a relatively long time to complete. 

The solution to this is to use the technique of simulating spiral waves in a frame of reference
that is moving with the tip of the wave. Foulkes and Biktashev~\cite{Foul09} published a method that could 
achieve not only simulations for RW but also MRW, something that other authors were not 
able to achieve \cite{Beyn-Thummler-2004, Hermann-Gottwald-2010}. 
This means that we could afford a much smaller numerical box in which to conduct the simulations
as the tip of the wave never reaches the boundaries of the box and, given a large enough box, the 
boundaries will not have an influence on the resulting spiral wave.

We review below the main results from Foulkes and Biktashev. Further details relating to these 
methods can be found in Foulkes and Biktashev \cite{Foul09}.

Let us consider the reaction-diffusion system of equations,
\begin{equation}
\frac{\partial {\bf{u}}}{\partial t}= \bf{D}\nabla^2 u+f(u),
\label{eqn:4a}
\end{equation}
where ${\bf{u}}={\bf{u}}({\bf{r}},t)=(u_1,u_2,\ldots,u_l)^\top\in\mathbb{R}^l$, $l\geq 2$, ${\bf{r}}=(x,y)^\top\in\mathbb{R}^2$,
${\bf{f}}\in{\mathbb{R}}^l$, and ${\bf{D}}\in{\mathbb{R}}^{l\times l}$ is the matrix of diffusion coefficients.
Foulkes and Biktashev considered a RDS that contained symmetry breaking perturbations,
which forced the spiral wave solution to drift. Since drift is of no concern in this work,
we consider the unperturbed RDS.

This system is invariant under the Euclidean group $SE(2)$, the group of the isometric transformations 
of the plane $\mathbb{R}^2\rightarrow\mathbb{R}^2$. This means that if ${\bf{u}}({\bf{r}},t)$ 
is a solution to equation (\ref{eqn:4a}), then $\tilde{{\bf{u}}}({\bf{r}},t)$ is another solution 
to (\ref{eqn:4a}) which is given by
\begin{equation*}
\tilde{{\bf{u}}}({\bf{r}},t)=T(g){\bf{u}}({\bf{r}},t),~~\forall g\in SE(2),
\end{equation*} 
where action $T(g)$ of $g\in SE(2)$ on the function ${\bf{u}}$ is defined as
\begin{equation*}
T(g){\bf{u}}({\bf{r}},t)={\bf{u}}(g^{-1}{\bf{r}},t).
\end{equation*}

It has been shown that elements of $SE(2)$ are concerned with the orientation-preserving transformations, 
i.e. only the spatial transformations, not the temporal transformations. The rigidly rotating spiral wave 
are independent of time in an appropriate rotating frame of reference but in this study, we are also 
interested in the meandering spiral wave solutions. So, for our current approach, the isotropy group of 
spiral wave solutions is trivial. This states that these solutions are invariant under any transformation.

By considering the spiral wave solution to (\ref{eqn:4a}) in an appropriate Banach space, and 
splitting out the motion of the spiral wave across and perpendicular to an appropriate representative manifold,
defined such that the tip of the spiral wave always remains on it, then we obtain the following system, 
\begin{eqnarray}
&&\frac{\partial\bf{v}}{\partial t} =\bf{D}\nabla^2\bf{v}+f(v)+(\bf{c},\nabla)\bf{v}+\omega\partial_\theta\bf{v},\label{eqn:4z}\\
&&v_{l_1}({\bf{0}},t)=u_*,~~v_{l_2}({\bf{0}},t)=v_*,\label{eqn:a1}\\
&&\frac{\partial v_{l_3}}{\partial x}({\bf{0}},t) = 0,~~\frac{\partial v_{l_3}}{\partial y}({\bf{0}},t) > 0,\label{eqn:a2}\\
&&\dv{\theta}{t} = \omega,~~\dv{{\bf{R}}}{t} = e^{\tau\theta}{\bf{c}},\label{eqn:a3}
\end{eqnarray}
where ${\bf{v}}={\bf{v}}({\bf{r}},t)=(v_1,v_2,\ldots,v_l)^T\in\mathbb{R}^l$ is the spiral wave solution 
in a frame of reference that is moving with the tip of the spiral wave, ${\bf{c}}(t)=(c_x(t),c_y(t))$ is the 
translational velocity of the spiral wave and $\omega(t)$ is its rotational velocity. The position and orientation of the tip
are given by ${\bf{R}}$ and $\theta$ respectively, hence equations (\ref{eqn:a3}) are equations of 
motion of the tip of the spiral wave. The fixed parameter, $\tau$, is the 
matrix $\tau=\begin{bmatrix}0 & -1\\1 & 0\end{bmatrix}$, meaning that exp($\tau\theta$) is the rotation matrix 
where the rotation is by angle $\theta$. 

We note that equation (\ref{eqn:4z}) is a 
reaction-diffusion-advection system of equations~\cite{Foul09, bik96}, whose spiral wave solutions ${\bf{v}}({\bf{r}},t)$ 
are such that their tip remains on the manifold. The conditions (\ref{eqn:a1}) \& (\ref{eqn:a2}) are the tip 
pinning conditions. The tip can be pinned at any point within the numerical box, but the definition here is such 
that they are pinned at the origin, which we place at the center of the numerical box, 
at an orientation determined by (\ref{eqn:a2}). 
\subsection{Reaction Kinetics}
We will be using the FitzHugh-Nagumo (FHN) model for the simulations within the work presented here.
This is a two variable, parameter dependent reaction-diffusion type model. Since this means that $l=2$, 
then we let $u_{l_1}=u_{l_3}=u_1$ and $u_{l_2}=u_2$.
In the equation (\ref{eqn:4a}), ${\bf{f}}({\bf{u}})$ defines the model kinetics, which, for the FHN model, are given by,
\begin{equation}
{\bf{f}}({\bf{u}}):\begin{bmatrix} u_1\\ u_2 \end{bmatrix}\mapsto\begin{bmatrix} \frac{1}{\epsilon}\left(u_1-\frac{u_1^3}{3}-u_2\right)\\ \epsilon\left(u_1+\beta-\gamma u_2\right)\end{bmatrix}.
\label{eqn:a7}
\end{equation}

We see that there are two variables - $u_1({\bf{r}},t)$, the excitation variable, and $u_2({\bf{r}},t)$, 
the inhibitory variable - together with three parameter, $\beta$, $\gamma$, and $\epsilon$.

The parameters are varied to give a variety of solutions. Winfree~\cite{Winfree-1991}
illustrated the spiral wave solutions in a parametric portrait, based on fixing the parameter 
$\gamma=0.5$, and varying the remaining two parameters, $\beta$ and $\epsilon$, to get a plethora of
solutions within a section of the parameter space, including regions of \emph{hypermeander} and plane waves. 
In general, we usually have $|\epsilon|\ll1$.
\subsection{Numerical Implementation}

Numerical implementation of this system is also detailed in~\cite{Foul09} and resulted in software called 
\textit{EZRide}~\cite{EZRide}. Operator splitting was utilised to simplify the otherwise complicated 
equations. We can rewrite equation (\ref{eqn:4z}) as,
\begin{equation*}
\frac{\partial {\bf{v}}}{\partial t}=\mathcal{F}[{\bf{v}}]+\mathcal{A}[{\bf{v}};{\bf{c}},\omega],
\end{equation*}
such that,
\begin{eqnarray*}
\mathcal{F}[{\bf{v}}] &=& \bf{D}\nabla^2 {\bf{v}}+{\bf{f}}({\bf{v}}),\\
\mathcal{A}[{\bf{v}};{\bf{c}},\omega] &=& ({\bf{c}},\nabla){\bf{v}}+\omega\frac{\partial {\bf{v}}}{\partial \theta},\\
                              &=& (c_x-\omega y)\frac{\partial {\bf{v}}}{\partial x}+(c_y+\omega x)\frac{\partial {\bf{v}}}{\partial y}.
\end{eqnarray*}					

For discretisation, we have a constant time step, $\Delta_t$, and space step, $\Delta_x$, covering the 
square spatial domain $(x,y)\in[-L/2,L/2]^2$, where $L$ is the length of the box in space units. 
The domain is divided into smaller squares by dividing the $x$ and $y$ axes
into $N_x$ and $N_y$ intervals respectively. For our purposes, we let
\[N_x = N_y = N = L/\Delta_x.\]
This means that there will be $N+1$ points along each axis. 

Let $\hat{\mathcal{F}}$ and $\hat{\mathcal{A}}$ be discretisations of $\mathcal{F}$ and $\mathcal{A}$ respectively. Our numerical computations are as follows,
\begin{eqnarray*}
\hat{V}^{k+\frac{1}{2}} &=& \hat{V}^{k}+\Delta_t\hat{\mathcal{F}}\left(\hat{V}^k\right),\\
\hat{V}^{k+1} &=& \hat{V}^{k+\frac{1}{2}}+\Delta_t\hat{\mathcal{A}}\left(\hat{V}^{k+\frac{1}{2}},\hat{\vec{c}}~^{k+1},\hat{\omega}^{k+1}\right),\\	
\hat{\theta}^{k+1} &=& \hat{\theta}^{k}+\Delta_t\omega^{k+1},\\
\hat{R}^{k+1} &=& \hat{R}^{k}+\Delta_t e^{\hat{\gamma}\hat{\theta}^{k+1}}\hat{\vec{c}}~^{k+1}.
\end{eqnarray*}
The timestep, $\Delta_t$, is given by
\begin{equation*}
\Delta_t = \frac{t_s\Delta_x^2}{4},
\end{equation*}
where $t_s$ is the ratio of the timestep to the diffusion stability limit, usually taken to be $t_s=0.1$ \cite{Foulkes-2009}.

\subsection{Reaction-diffusion step}
Foulkes et al~\cite{Foul09} and Barkley~\cite{Bark91} used the initial steps for computation just as the same as used in the Barkley's EZ-SPIRAL software. 
Further, they added more numerical computational steps to it. So, for the reaction diffusion part, a first order accurate forward Euler method was used 
to calculate the temporal derivatives, and the Five Point Laplacian method for the Laplacian.
\subsection{Advection step; to calculate $c_x$, $c_y$ and $\omega$}
An upwind second-order accurate approximation of the spatial derivatives is used in $\hat{\mathcal{A}}$. In this step, the discretisation of $\hat{V}^{k+1}$ 
at the tip pinning points is used to calculate $\hat{\vec{c}}_x~^{k+1}$, $\hat{\vec{c}}_y~^{k+1}$ and $\hat{\omega}~^{k+1}$ so that after every step, the tip 
pinning conditions are correctly satisfied.
\subsection{Tip pinning conditions}
Pinning the tip of a meandering spiral wave was achieved by choosing two isolines, one for excitation and one for inhibitory, whose values are located within 
the range of values of both excitation and inhibitory variables. The values can be determined by considering the phase portrait relating to the kinetics
used. The full details of the choice of tip pinning conditions are given in Foulkes \& Biktashev~\cite{Foul09}, and the reader should refer these full 
explanations. 

In summary, discretisation of the tip pinning conditions lead to
\begin{eqnarray*}
v_{l_1}^{(i_0,j_0);k} &=& u_*\\
v_{l_2}^{(i_0,j_0);k} &=& v_*\\
v_{l_1}^{(i_0+i_{inc},j_0+j_{inc});k} &=& u_*\\
v_{l_1}^{(i_0+i_{inc},j_0+j_{inc});k} &<& v_*
\end{eqnarray*}
where the subscripts are the variable identifier, and the super scripts represent the spatial and temporal variables respectively. 
The choice of tip pinning condition ensures that the tip of the spiral wave is fixed at a certain orientation and position regardless
of whether the spiral wave is rigidly rotating or meandering. 
\subsection{Reconstruction of tip trajectory}
The EZRide software has an in-built algorithm for reconstructing a tip in the comoving FOR~\cite{Foul09}. Further, to check our calculations and draw the tip trajectories, we considered equation (\ref{eqn:a3}) and solved it using the numerical scheme which is given as
\begin{eqnarray}
\theta_{k+1} &=& \theta_{k}+\Delta t \omega_{k}\label{eqn:r01},\\
x_{k+1} &=& x_{k}+{\Delta t}((c_x)_k\cos\theta_k-(c_y)_k\sin\theta_k)\label{eqn:r02},\\
y_{k+1} &=& y_{k}+{\Delta t}((c_x)_k\sin\theta_k+(c_y)_k\cos\theta_k)\label{eqn:r03}.
\end{eqnarray}

Therefore, the equations (\ref{eqn:r01}), (\ref{eqn:r02}) and (\ref{eqn:r03}) represents the reconstructed tip trajectory. 

\subsection{Other details}
Throughout the studies conducted here, we used Neumann boundary conditions. Several authors have noted the effect 
of the boundaries upon the behaviour of the spiral waves, and different boundary conditions can lead to different
solutions. However, it has been noted~\cite{bik03,bika09,Foul10} that provided the tip of the spiral is sufficiently 
far from the boundary, the boundary will not affect the overall dynamics of the spiral wave determined by its
tip. This in turn is related to the response functions of the spiral wave, which are localised at the tip \cite{bik03}. 
Therefore, we can use either Neumann or Dirichlet boundary conditions, or a mixture of both, within our simulations. 

As a guide to what sufficiently far from the boundary means, we usually take this to mean that there is a full
wavelength between the tip and the boundary, i.e. the distance measured from the tip of the wave to the part of the 
arm of spiral that has the same $u_1$ and $u_2$ values as at the tip and has rotated around by $2\pi$.

Another consideration of the numerical implementation of the system~(\ref{eqn:4z}-\ref{eqn:a2}, \ref{eqn:a7}) is that
of stability. It was shown~\cite{Foul09} via a von Neumann analysis that we require,
\[
|c_x| < \frac{\Delta^2_x}{2\Delta_t},\ 
|c_y| < \frac{\Delta^2_y}{2\Delta_t},\ 
|\omega| < \frac{1}{N_x\Delta_t},
\]
for stability in the calculation of $c_x$, $c_y$ and $\omega$.
If the values of $c_x$ and $c_y$ went beyond these limits, then their values were restricted to the limits. 
One of the techniques employed in order to reduce stability was that as the advection terms are ``switched on'',
then the tip of the spiral was moved immediately to the tip fixation point. It is also important to note that here $\omega$ was not calculated initially
until the orientation of the spiral wave was met, at which point $\omega$ was given an initial value of zero, and 
the calculated as usual. Foulkes \& Biktashev discovered that $\omega$ is very sensitive to calculations, which
led to the definition of the tip pinning conditions. Unlike $c_x$ and $c_y$, if $\omega$ went beyond its limiting
value, then it was allocated a value of zero rather than its limiting value. Allocating $\omega$ its limiting
value eventually led to instabilities within the system.

\section{Numerical bifurcation approach}
We aim to show the nature of the bifurcation responsible for the transition of spiral waves from rigid rotation
to meander by generating solutions in a frame of reference that is comoving with the tip of the spiral wave. 
This means that even for large core spiral waves, we can still afford a relatively small computational space. 
Furthermore, Foulkes \& Biktashev~\cite{Foul09} showed that within the solutions to~(\ref{eqn:4z}-\ref{eqn:a2}, \ref{eqn:a7}),
the \emph{quotient data}, consisting of the dynamic variables ${\bf{c}},\ \omega$, form limit cycle solutions.
We therefore study the growth of these limit cycle solution from the onset of meander,
and the nature of the growth of this data will indicate the type of bifurcation taking place.
\subsection{Methodology: general overview}
\begin{figure}[tb]
\begin{center}
\centering
\includegraphics[width=0.65\textwidth,trim={8.5cm 0cm 0cm 0cm},clip]{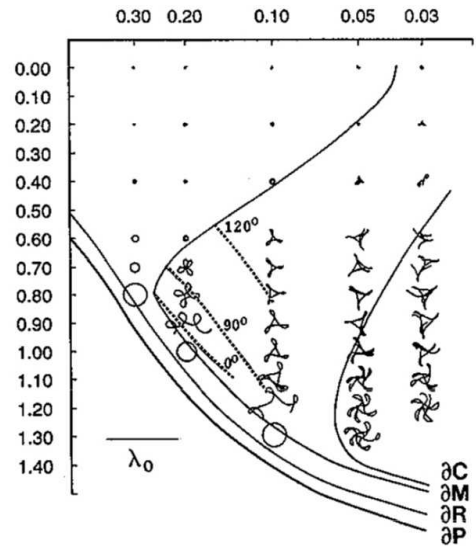}
\caption{Parametric Portrait for FHN for $\gamma=0.5$~\cite{Winfree-1991}.}
\label{fig:fhn}
\end{center}
\end{figure}

Consider the parameter portrait from Winfree for FHN, as shown in figure (\ref{fig:fhn}).
We see that there are different regimes of types of motion of spiral waves, according to values of $\beta$ and $\epsilon$.
We decided to study the growth of the limit cycles relating to meandering spiral waves, by analysing the quotient
data for a range of spirals which, on varying one of the parameters, go from rigid rotation to meander and then back to rigid rotation.
From figure (\ref{fig:fhn}), we decided to fix $\epsilon=0.2$, and starting at $\beta=0.570$ within the upper rigidly rotation
space, we varied the $\beta$-parameter in steps of $\Delta_\beta=0.001$ to get a spiral wave solution for each value of $\beta$. 
Our initial thoughts were that this choice of $\Delta_\beta$ was sufficient to generate a range of 
solutions which show the nature of the bifurcation. As we will see in later sections, the choice of $\Delta_\beta$ will lead to qualitatively
similar results, but quantitatively different ones, showing sensitive dependence on initial conditions. 

\begin{figure}[tb]
\begin{center}
\begin{minipage}[htbp]{0.48\linewidth}
{\centering
\includegraphics[width=0.95\textwidth]{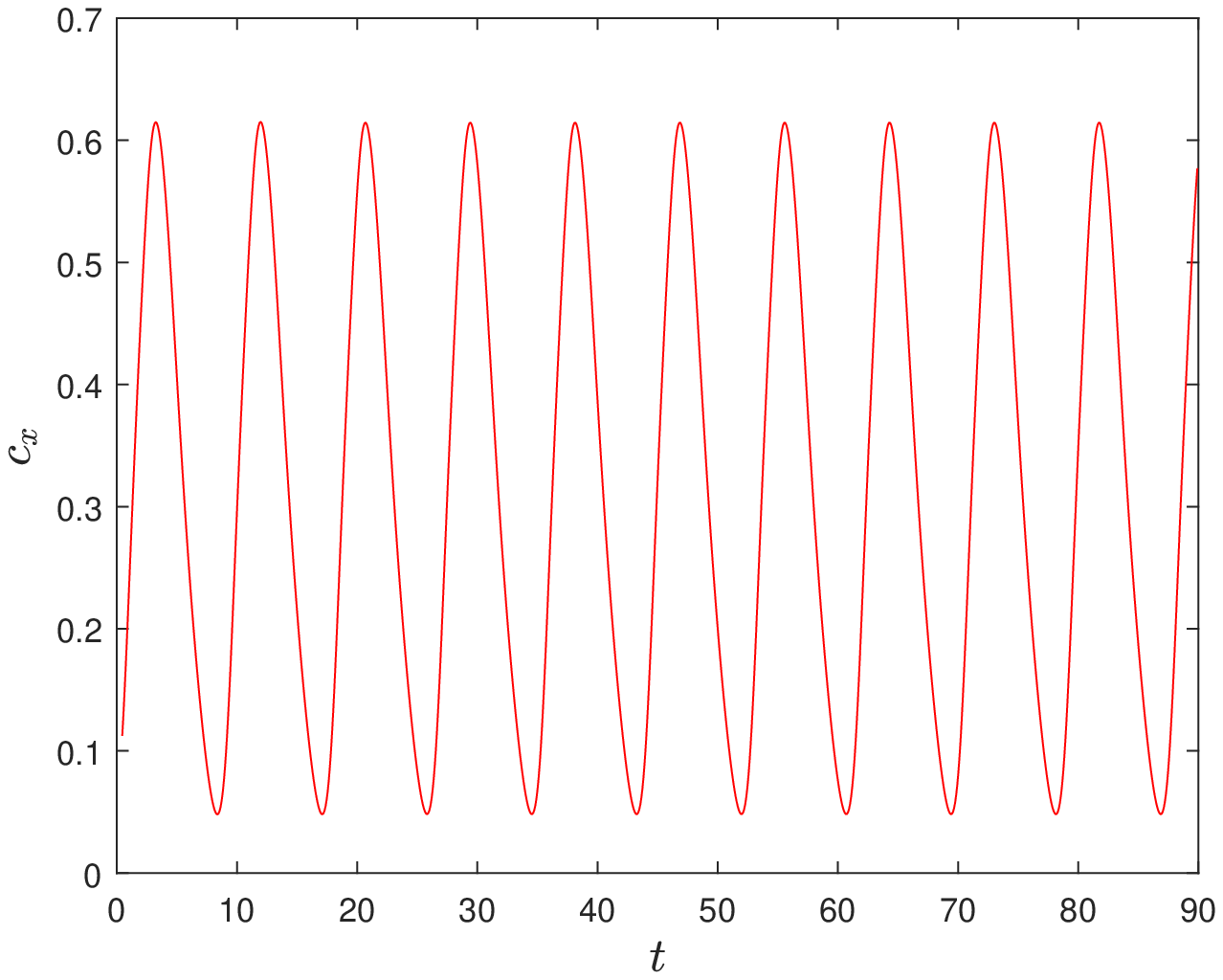}}
\end{minipage}
\begin{minipage}[htbp]{0.48\linewidth}
{\centering
\includegraphics[width=0.95\textwidth]{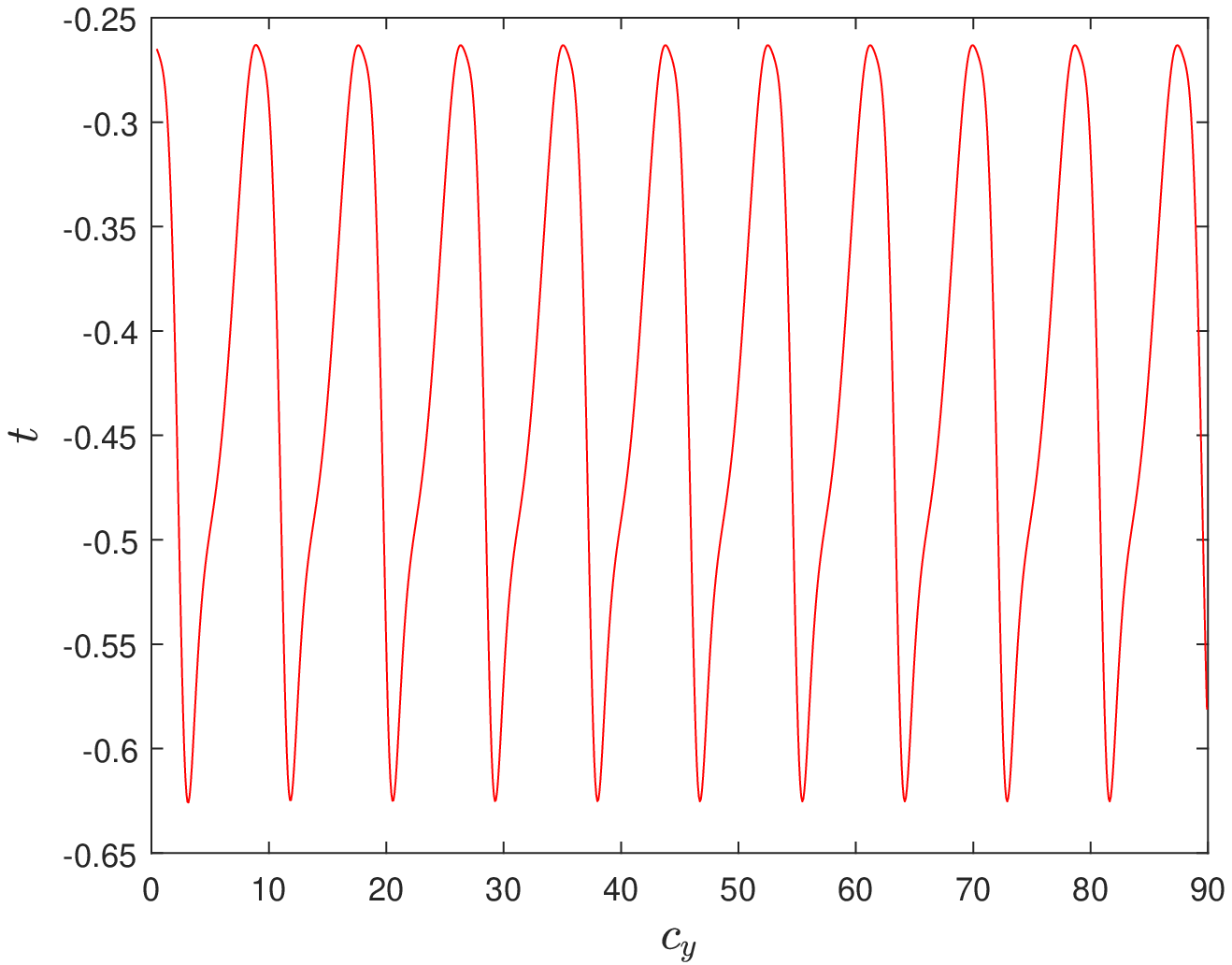}}
\end{minipage}
\begin{minipage}[htbp]{0.48\linewidth}
{\centering
\includegraphics[width=0.95\textwidth]{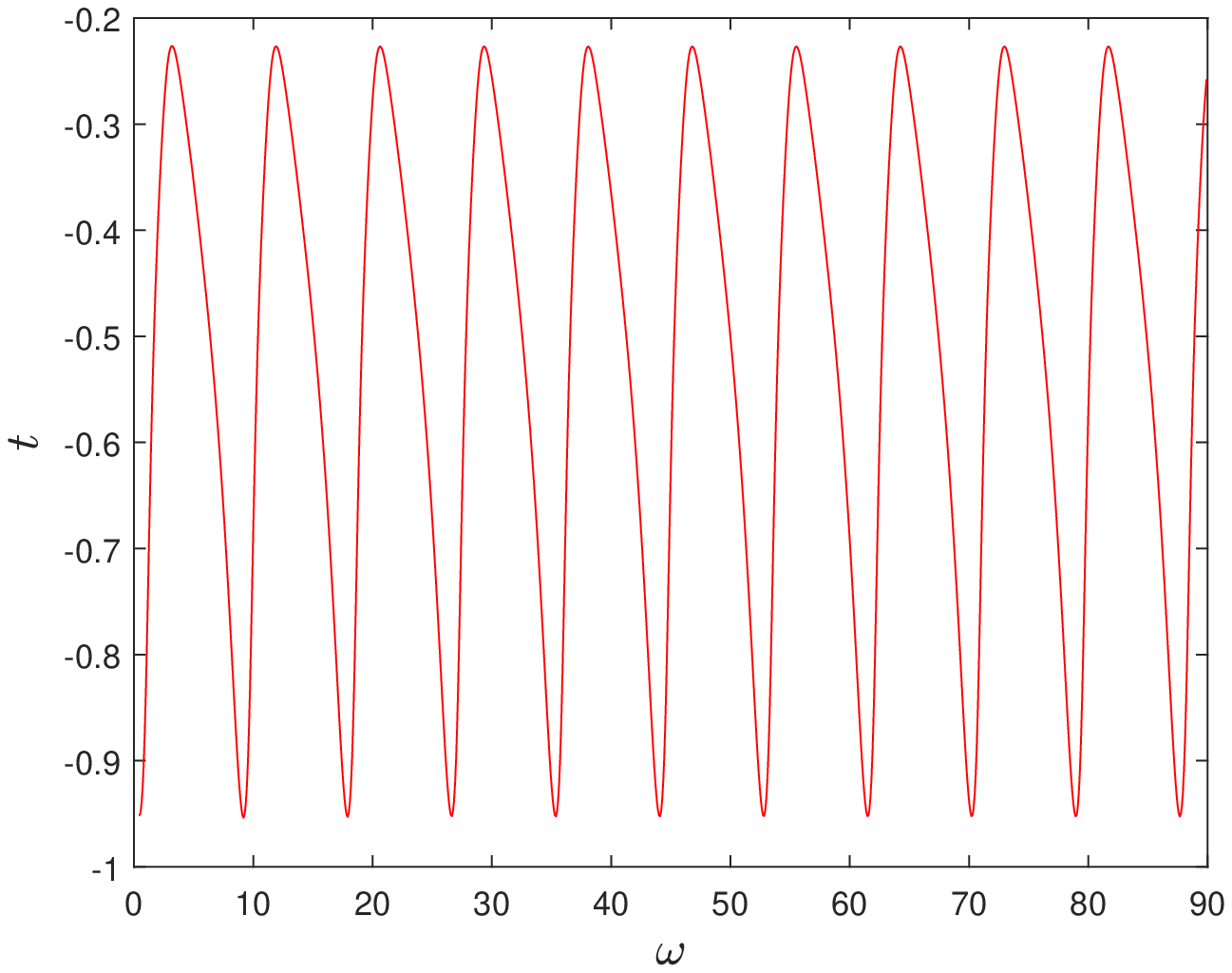}}
\end{minipage}
\begin{minipage}[htbp]{0.48\linewidth}
{\centering
\includegraphics[width=0.95\textwidth]{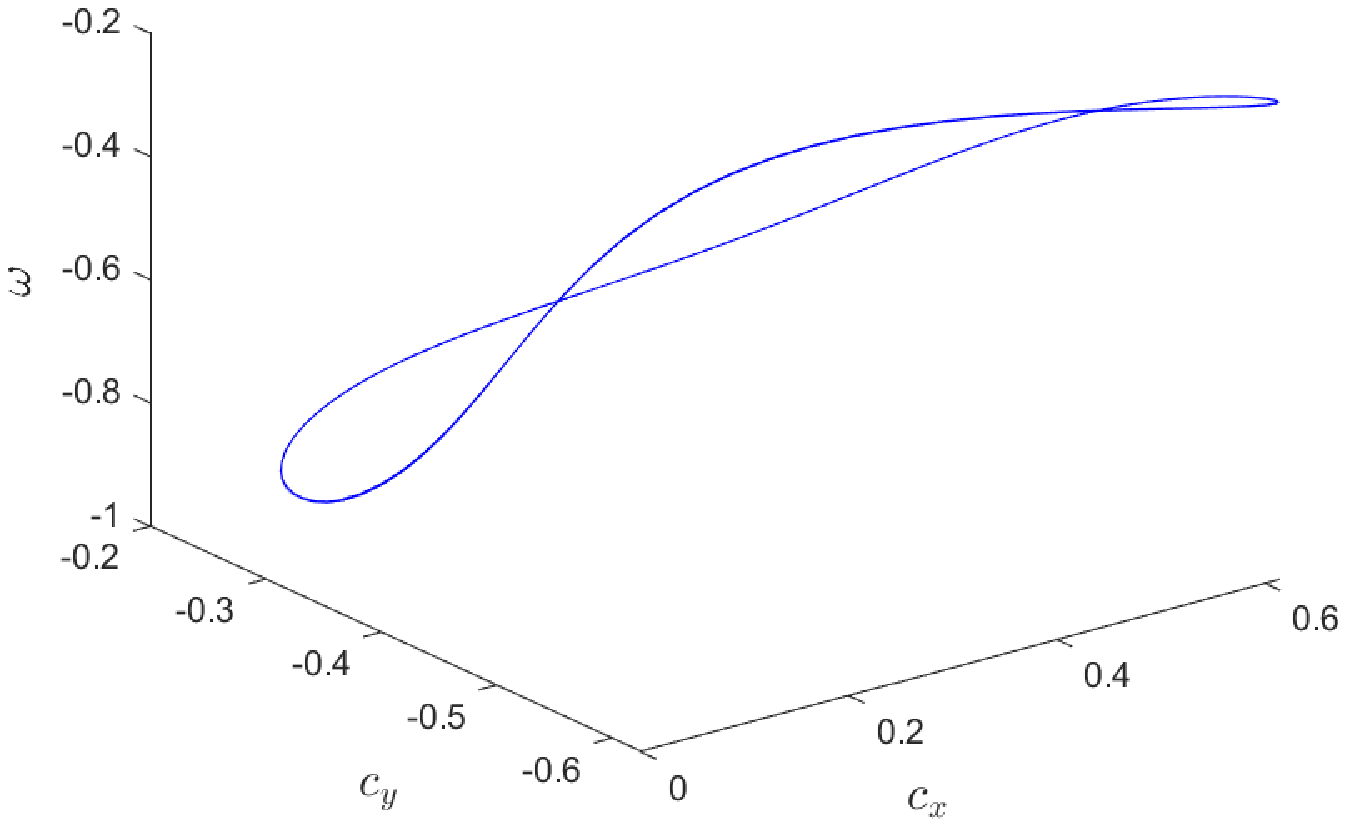}}
\end{minipage}
\caption{Quotient data for a meandering spiral wave with $\beta=0.751$. 
Time v/s (top left) $c_x$, (top right) $c_y$, (bottom left) $\omega$. The full
limit cycle is shown (bottom right).} 
\label{fig:graphs_fig}
\end{center}
\end{figure}

Each simulation records not only the quotient data, but also the final conditions
from the end of the simulation. The initial conditions for each new simulation is taken 
as the same as the final conditions of the previous simulation.
We shall see that, when analysing the results, we need to take into account any transient part of the current solution when 
using the final conditions of the previous solution as initial conditions for the current solution.
Although this transient period can sometimes vary in length, only the data taken after at least five complete periods had
occurred. In some, but not many, cases, it was obvious that the transient period occurred for more that five periods and 
therefore the data for those simulations was analysed to see where the transient period had ceased. The remaining data
was then analysed in the usual way. 

Once all the simulations had been completed, the quotient data was then analysed and the ``size'' of the limit cycles were then plotted
against the parameter, $\beta$, which we shall also call the \emph{bifurcation parameter}. To do this, we need to define this size and 
so introduce below the \emph{quotient size}, $Q_s$, of the limit cycle of a spiral wave.

Furthermore, we shall consider the bifurcation points which are the values of the $\beta$ parameter at which the bifurcation takes place,
and relate directly to the last simulation for which there is no limit cycle.

\subsection{Quotient size}

The shape of the limit cycle in figure (\ref{fig:graphs_fig}) is irregular. In previous studies~\cite{Bark94a}, the radius of the limit
cycles were measured, but this technique cannot be applied here, due to the irregularity of the shape of the limit cycle. We therefore
calculated the ``distance'' around the limit cycle from one point on the cycle all the way around back to that point by calculating
the arc length.

We know that arc length of a function ${\bf{q}}(t)$ in the interval $t\in[t_1,t_2]$ is given by,
\[Q_s = \int_{t_1}^{t_2}||\dot{{\bf{q}}}(t)||\, dt.\]
For the spiral wave limit cycles, we let ${\bf{q}}(t)=\left(c_x(t),c_y(t),\omega(t)\right)$, and if we take
the integral over one whole period of the spiral wave, $T$, then 
\begin{equation*}
Q_s = \int_{t}^{t+T}\sqrt{\dot{c}^2_x+\dot{c}^2_y+\dot{\omega}^2}\, dt.
\end{equation*}
Since we do not know the exact form of $c_x$, $c_y$ and $\omega$, then we need to use the discretised form of the arc length
formula,
\begin{eqnarray}
Q_s &\approx& \sum_{i=j}^{j+N}\sqrt{(\hat{c}^{i+1}_x-\hat{c}^i_x)^2+(\hat{c}^{i+1}_y-\hat{c}^i_y)^2+(\hat{\omega}^{i+1}-\hat{\omega}^i)^2},\nonumber\\
&& 
\label{eqn:qs}
\end{eqnarray}
where $N=T/\Delta_t$, $j$ is a starting point on the limit cycle, and $\hat{c}^i_x$ is the discretised value of $c_x$ at the $i$-th step, 
with similar notations for $c_y$ and $\omega$..

As noted earlier, due to the transient period of the spiral wave when the simulation first starts, we neglect the first five
periods of the simulation. If there are $n$ full periods left of the simulation, the quotient size of 
those $n$ periods were calculated separately 
using equation (\ref{eqn:qs}), and then $Q_s$ will be the average of those $n$ quotient sizes. 
 
In the case where there is rigid rotation, then $Q_s=0$ since once the transient period has passed, 
the values of $c_x(t)$, $c_y(t)$, and $\omega(t)$
all remain constant. Constant quotient is indeed an indication of rigid rotation, and once this constancy is detected, then the next simulation is
started after 50 time steps. This ensured that the simulation had settled to a solution that was rigidly rotating.

\section{Results}
For all the initial simulations, we have used $L=30$, $\Delta _x=\frac{1}{5}$, $\Delta_t=0.001$ and $t_s=0.1$. 
These numerical parameter values were carefully taken from Foulkes et al, 2009, so that we get not only accurate simulations
but computational fast generation of these simulations \cite{Foul09}. In the comoving FOR, we 
observed that for rigidly rotating spiral waves, our solution becomes stationary and that the quotient data stabilises to
constant values for $c_x$, $c_y$, and $\omega$. As $\beta$ varies, quotient system is no longer constant and 
is in fact periodic. This corresponds to meandering spiral wave solutions, exhibiting complicated quasiperiodic motion. 

\begin{figure}[tbph]
\begin{center}
\centering
\includegraphics[width=0.53\textwidth,trim=0mm 18mm 0mm 0mm,clip]{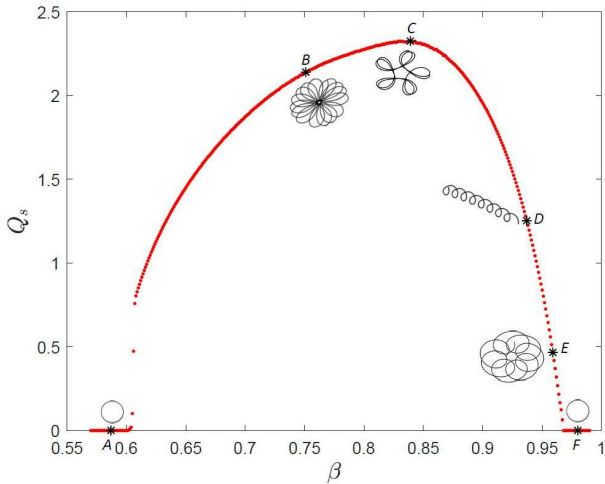}
\caption{The $\beta$-$Q_s$ plot. Each dot represents a single single simulation. Shown are 
reconstructed tip trajected from the quotient data. Simulations $A$ and $F$ are RW; $B$ and $C$ 
are MRW with outward facing ``petals''; $D$ is near the 1:1 resonance line and has an extremely
large core radius (something which previous authors could not simulate quickly); and $E$ is 
MRW with inward facing ``petals''.}
\label{fig:Arc_Length}
\end{center}
\end{figure}

The initial results are shown in figure (\ref{fig:Arc_Length}), where the parameter, $\beta$, is plotted against the 
quotient size, $Q_s$. The simulations started from $\beta=0.570$ and $\beta$ increased with
a step size $\Delta_\beta=0.001$. The quotient size was zero for $0.570\leq\beta\leq0.601$, 
indicating rigid rotation. For $0.602\leq\beta\leq0.967$, we had $Q_s>0$, meaning that the limit cycles were 
present and solutions are classed as meandering spiral waves. Increasing $\beta$ further, we found that for 
$0.968\leq\beta\leq0.990$, $Q_s=0$ meaning rigidly rotating spiral waves were present.

The change in the dynamic behaviour from RW to MRW is due to the Hopf bifurcations arising from the stable steady states, 
RW \cite{Bark94a}. These Hopf bifurcation take place at $\beta=0.601$ and $\beta=0.968$. 
Furthermore, as we increase $\beta$ from the critical point $\beta=0.601$, there is a 
sudden change in the qualitative behaviour of the system from a stable 
steady state (RW) to an oscillatory state (MRW).  

It is clear that as we vary $\beta$ within the MRW region, i.e. for $0.602\leq\beta\leq0.967$, then the 
growth of the limit cycles, as measured by $Q_s$, initially increases from zero up to a maximum, and then decreases again to 
zero. In order to see whether the arc formed by plotting $\beta$ against $Q_s$ contained any special features, 
such as the maximum of the arc relating to a specific type of meander, we chose some values of $\beta$ and
plotted their reconstructed trajectories from the quotient data. We illustrate the various types of 
solutions in figure (\ref{fig:Arc_Length}). There were no particular surprises in this analysis and
the results tied in with what Winfree observed in his parametric portrait of FHN that we reproduce in 
figure (\ref{fig:fhn}).

A feature of this plot is that as we increase $\beta$ from the left hand Hopf bifurcation point, there is a significant
gap between values of $Q_s$ which is not in line with the rest of the plot. We see that once $\beta$ has 
increased beyond 0.601, $Q_s$ grows slowly and then growth suddenly accelerates very quickly, before decelerating near
the peak of the plot. It then accelerates (this time in the negative direction of $Q_s$) down to the 
right hand bifurcation point. We therefore decided to look closer at this gap near the left hand bifurcation
point.

According to the previous studies which were conducted using Barkley's model~\cite{Bark94a} or the Belousov-Zhabotinksy 
reaction~\cite{Swinney}, it was observed that a supercritical Hopf bifurcation is responsible for the 
transition from rigid rotation to meander.  
At the left hand bifurcation point, we observe a discontinuous jump in the growth of quotient size. 
We can also see that immediately after the bifurcation, there is no square root relationship between $\beta$
and $Q_s$. A square root was attempted to be fitted to the arc immediately after the bifurcation point, 
but it would not fit in perfectly. Thus, the initial observation is that result does not tie with 
the analysis of supercritical Hopf bifurcations due to the absence of square root characteristic. 

However, the discontinuous jump observed in the bifurcation diagram depicted in the figure (\ref{fig:Arc_Length}), 
signals a subcritical Hopf bifurcation~\cite{izhikevich}. It is also known that if we vary our bifurcation 
parameter back and forth across the Hopf bifurcation point, we wouldn't expect to jump back to the same value of 
$\beta$, where it lost its stability. We therefore decided to run the simulations again, but this time starting
at $\beta = 0.990$. If there is a subcritical Hopf bifurcation present, then we should observe hysteresis,
which is associated with the bistable region~\cite{izhikevich}. 

\section{Hysteresis}\seclabel{hyst}

\begin{figure}[tbp]
\begin{center}
\centering
\includegraphics[width=0.5\textwidth]{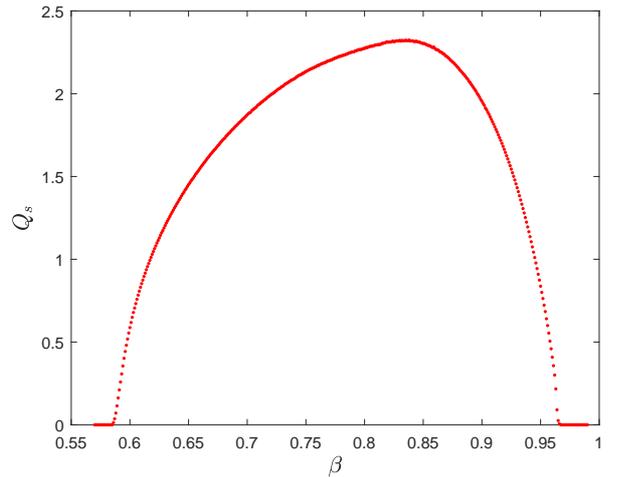}
\caption{Bifurcation diagram for the model parameter $\beta$ for the reverse run across the $\beta$-range from $0.990$ to $0.570$ with the $\beta$-step of $0.001$.}
\label{fig:Arc_back}
\end{center}
\end{figure}

\begin{figure*}[t]
\begin{center}
\begin{minipage}{0.49\linewidth}
\centering
\includegraphics[width=0.9\textwidth]{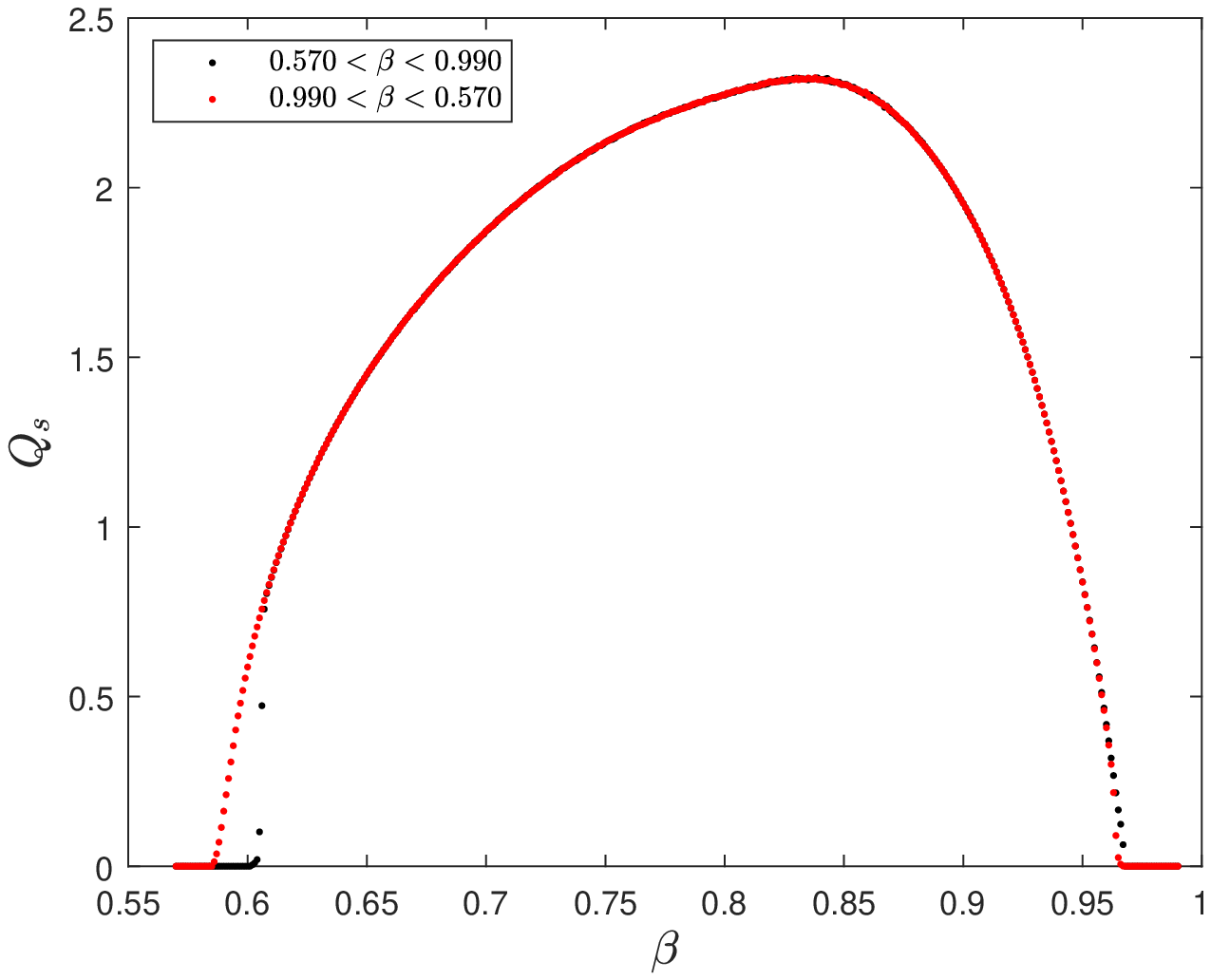}
\end{minipage}
\begin{minipage}{0.49\linewidth}
\centering
\includegraphics[width=0.9\textwidth]{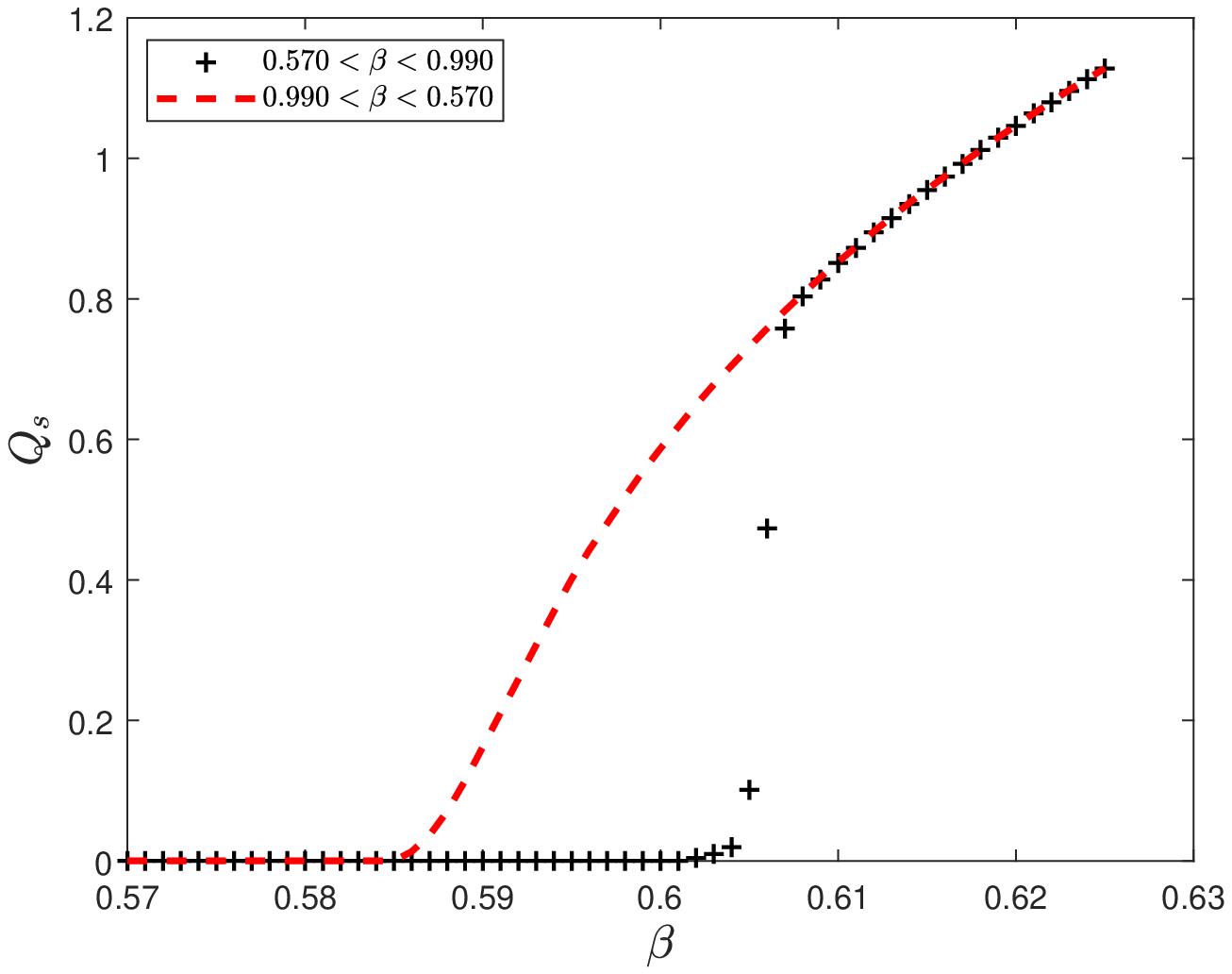}
\end{minipage}
\caption{Bifurcation diagram: the $\beta$-$Q_s$ plots depicting the hysteresis region. The {\textcolor{black} {\bf black}} curve represents forward run whereas {\textcolor{red} {\bf red}} curve represents reverse run across the chosen $\beta$-range. Both 
figures are from the same set of data with the full set of data shown (left) and the data around the hysteresis region
shown (right). Dots were used for (left) to illustrate the coinciding of the solutions for most of the data, and 
{\textcolor{red} {\bf red}} lines and {\textcolor{black} {\bf black}} crosses for the smaller data region (right).}
\label{fig:hyst}
\end{center}
\end{figure*}

The simulations were now run backwards from $\beta$ starting from $0.990$ and decreasing in step of 
$\Delta_\beta=0.001$ again to $0.570$. We performed the same calculations as in the previous section to 
calculate $Q_s$ for all the values of $\beta$ and observed its growth against $\beta$, as shown in figure 
(\ref{fig:Arc_back}). We note that the bifurcation diagram for the simulations run across backward looks 
similar to the one in figure (\ref{fig:Arc_Length}). However, for these simulations, the Hopf 
bifurcation point on the left-hand side in the reverse case has now shifted to $\beta=0.586$, whereas 
on the right-hand side the Hopf bifurcation point shifts very slightly to left at $\beta=0.967$.

Further, combining both the bifurcation diagrams figures (\ref{fig:Arc_Length}) and (\ref{fig:Arc_back}), we observed 
that both the graphs coincide coincide exactly except for only a few values of $\beta$, depicted in figure (\ref{fig:hyst}). 
This depicts the presence of the hysteresis region in which steady and oscillating states coexist. Hence, this region 
is associated with bistability. The presence of the hysteresis zone is an important characteristic for subcritical 
Hopf bifurcations, where the system can be in more than one state~\cite{gopalakrishnan, izhikevich, mergell}. 
In our study, the hysteresis region exhibits both RW and MRW solutions for a small range of $\beta$. Therefore, 
the existence of both the solutions and change in the $\beta$-value corresponding to the equilibrium clearly 
corresponds to the case of subcritical Hopf bifurcation.

\begin{figure*}[t]
\begin{center}
\begin{minipage}{0.49\linewidth}
\centering
\includegraphics[width=0.9\textwidth]{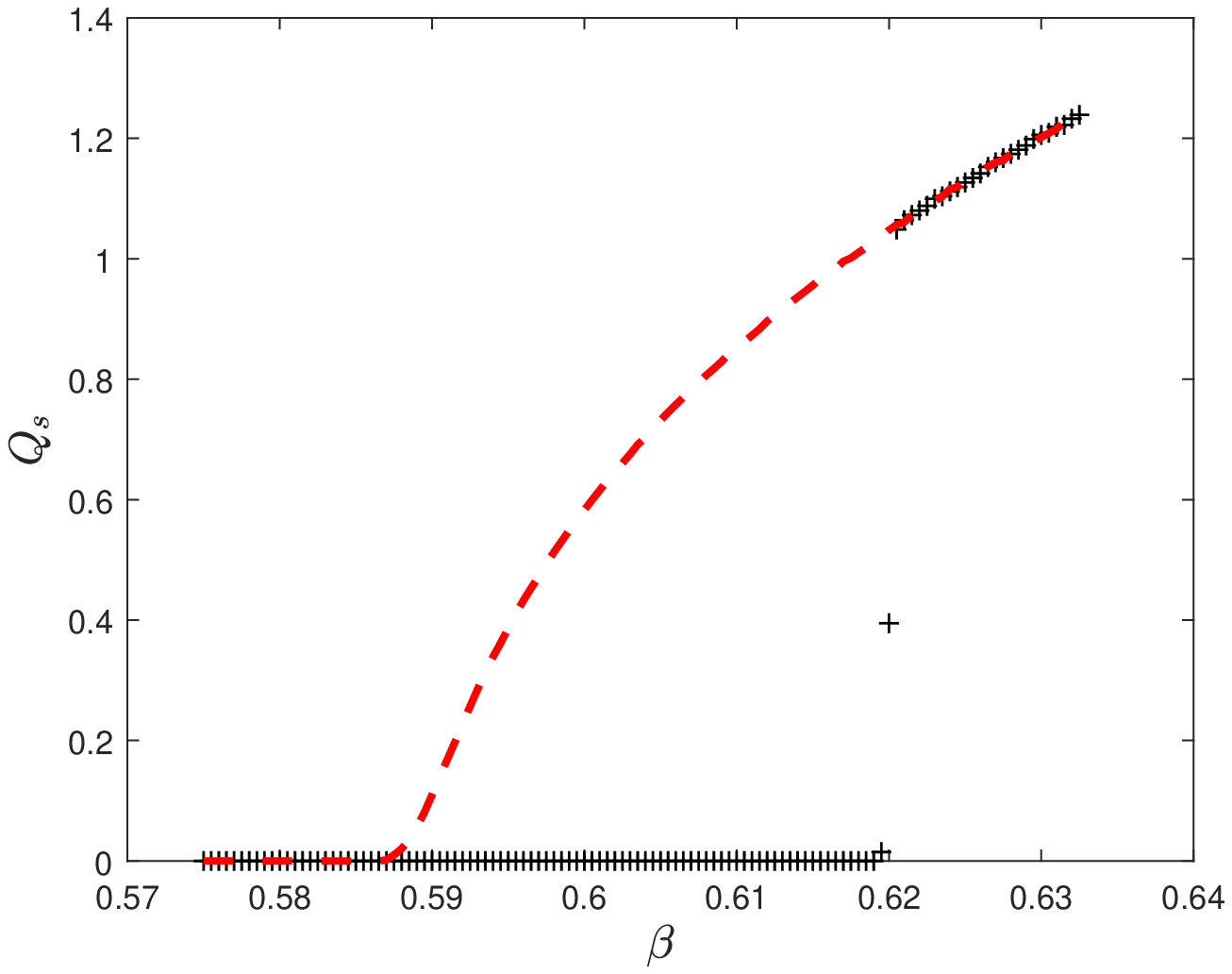}
\end{minipage}
\begin{minipage}{0.49\linewidth}
\centering
\includegraphics[width=0.9\textwidth]{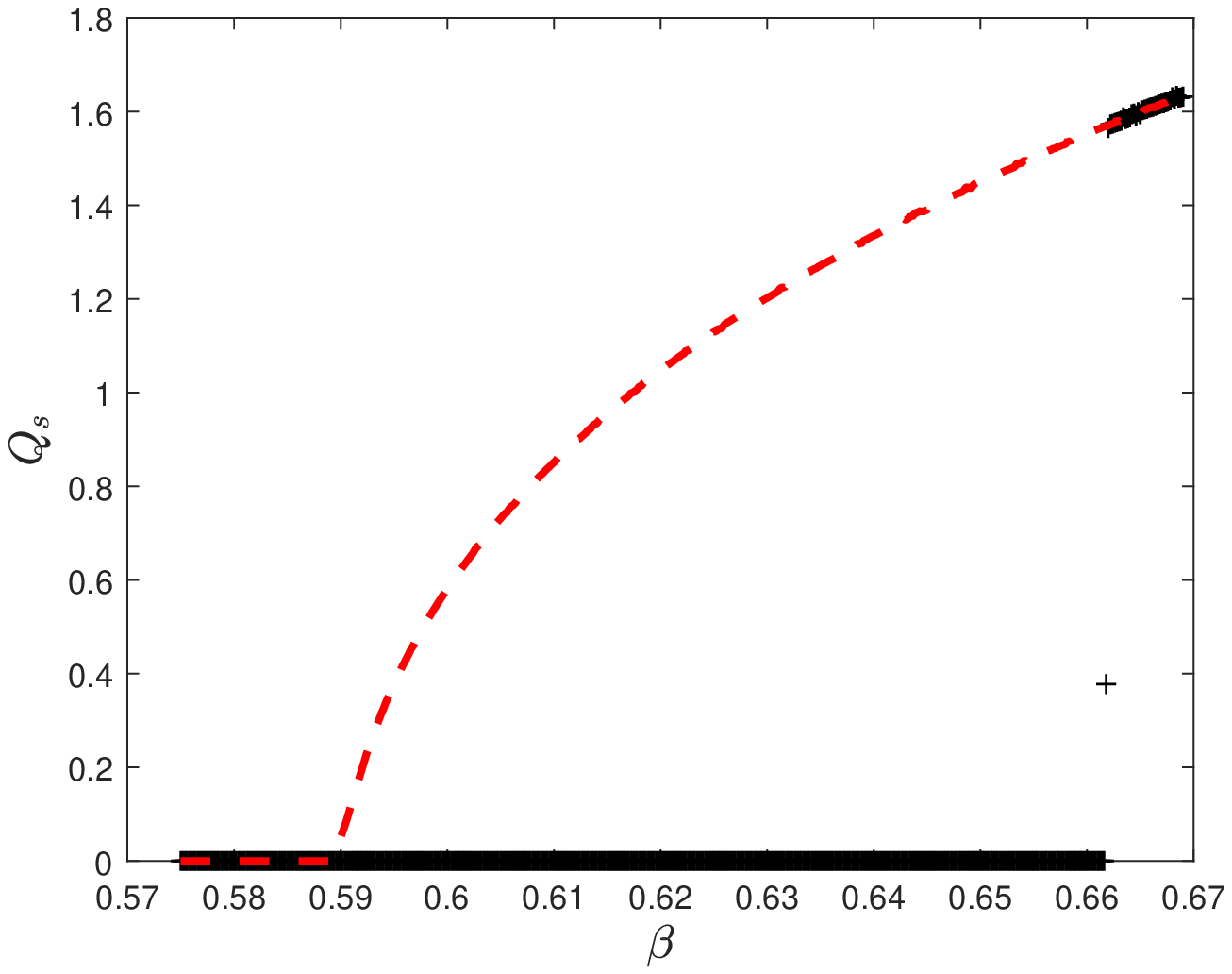}
\end{minipage}
\caption{Variation in $\Delta_\beta$: (left) 0.0005, (right) 0.0002.}
\label{fig:vary_b}
\end{center}
\end{figure*}

Now, concentrating only on the left hand side of the hysteresis region, we varied the step size of $\beta$. 
The idea behind this was to see if we generate further data points for the part of the hysteresis in which there is a jump, 
i.e. for when $0.604\leq\beta\leq0.607$.
All the other model and numerical parameters were kept same and we tested the hysteresis region with the smaller 
$\Delta_\beta$-steps. We decreased it to half the original step size of $\beta$ for one set of simulations
and further to one-fifth of the original $\Delta_\beta$-step for another set of simulations. Following the same 
procedure of calculating $Q_s$ and plotting it against $\beta$, it can be seen that the bifurcation points for 
the forward runs shift further right with the smaller $\beta$-steps, therefore increasing the width of the hysteresis 
zone. The growth of $Q_s$ against $\beta$ can be seen in figure (\ref{fig:vary_b}). 

One such explanation for this is that for larger $\beta$-steps, the step is large enough to induce a perturbation 
such that there is sensitivity to initial conditions causing a premature switch from RW to MRW motion. We note that the
procedure of conducting these simulations was to use the final conditions of the previous simulation as the initial conditions
of the next simulation. Thus, bistability within the system kicks in early when the $\beta$-step is too large. However, 
for smaller $\beta$-steps such as 0.0005 and 0.0002, the initial conditions for one simulation may be close to the final 
conditions of the previous previous such that the transition from RW to MRW is not influenced by the distance between 
the $\beta$-values. This describes the effect of sensitivity to initial conditions~\cite{feldman}. Furthermore, 
the bifurcation points for the forward and the reverse run across the $\beta$ values are highlighted in table (\ref{table:b_step}).

\begin{table}[tbph]
\begin{center}
\begin{tabular}{|m{2.5cm}|m{2.5cm}|m{2.5cm}|c|c|c|}
\hline
$\beta$-step & Forward run & Reverse run\\
\hline
0.001 & $0.602$ & $0.586$ \\
\hline
0.0005 & $0.6195$ & $0.587$ \\
\hline
0.0002 & $0.6618$ & $0.589$ \\
\hline
\end{tabular}
\caption{Variation in $\beta$-steps and their corresponding bifurcation points}
\label{table:b_step}
\end{center}
\end{table}

Further, analysing the behaviour across the hysteresis region, we decided to run the simulations with the same 
range of $\beta$ but now for a different value of epsilon. We chose $\epsilon=0.25$ and varied $\beta$ both forward 
and backward across the range. 

\begin{figure}[tbp]
\begin{center}
\centering
\includegraphics[width=0.53\textwidth]{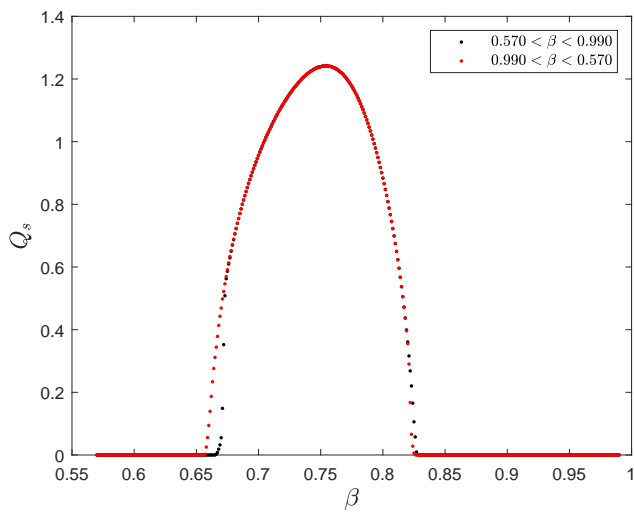}
\caption{Bifurcation diagram $\epsilon=0.25$ with $\beta$ varied across forward and backward over the range.}
\label{fig:Arc_eps}
\end{center}
\end{figure}

From figure (\ref{fig:Arc_eps}), we can see that both the graphs coincide with each other apart for the few values of 
$\beta$ near the bifurcation point, which seems more obvious on the left hand side, similar to our original result. 
The width of the hysteresis region represents the parameter range where the system is bistable. It can be clearly 
seen that the width of the hysteresis zone is also affected with the change in the value of $\epsilon$. Here, we 
discovered that the width of the hysteresis zone decreases with the increase in $\epsilon$. The presence of hysteresis 
also confirms that the transition to instability is subcritical. 

\section{Conversion}\seclabel{shock}

Our next concern was to consider only the hysteresis region and to check if it is possible, for a pair of solutions
corresponding to the same $\beta$ value, to convert from one solution to another. We only tested for the region 
on one side of the original bifurcation diagram ($\Delta_\beta$=0.001), being more prominent on the left hand 
side, shown in figure (\ref{fig:hyst}). 
We described the conversion of one type of solution to another under perturbation such as single shock 
defibrillation~\cite{Foul10, holdendefib, morgan2009}.

The main focus here is on the bistable region where there exists both steady and periodic states. In our case, we 
investigate a parameter region in an excitable media, where RW and MRW solutions coexist for the same parameter
values. They differ significantly as rigid motion depicts steady states whereas meander depicts oscillatory states.

It can be seen that for a particular value of $\beta$ within the hysteresis region given by $0.586<\beta<0.602$, 
depicted in the figure (\ref{fig:hyst}), we have two types of solutions. Here, we consider the 
transition between the two solutions and chose $\beta=0.595$, which exhibits steady state (RW) while we run the 
simulations forward across the $\beta$-range whereas it shows the periodic behaviour (MRW) with the reverse run.

Now, to test if the spiral wave solutions can be converted from one type to another, there are many ways in which
this can be achieved. We decided to use the method of single shock defibrillation by adding a perturbation to 
the system. Foulkes et al.~\cite{Foul10} have shown that we can convert one type of spiral into another 
by means of a single shock, and so we utilise techniques from their work.

We added a perturbation to the $u$-field, uniformly throughout space at a specified time 
$T$~\cite{Foul10}.  
\begin{equation*}
\frac{\partial{\bf u}}{\partial t} = {\bf D}\nabla^2{\bf v}+{\bf f}({\bf v})+{\bf h}(t),
\end{equation*}
where ${\bf u}=(u,v)$ and ${\bf h}(t)=(h_u(t), h_v(t))^\top$ is a time-dependent perturbation, defined as
\begin{equation*}
h_u(t) = A\delta(t-T),~h_v(t) = 0,
\end{equation*}
where $A$ is a constant.

Since we consider solutions in a comoving frame of reference, we consider solutions to the 
reaction-diffusion-advection system of equations~\cite{Foul09}
\begin{equation*}
\frac{\partial{\bf v}}{\partial t} = {\bf D}\nabla^2{\bf v}+{\bf f}({\bf v})+ \tilde{{\bf h}}(t)+({\bf c},\nabla){\bf v}+\omega\frac{\partial{\bf v}}{\partial \theta},
\end{equation*}
where $\tilde{{\bf h}}(t)$ is defined as
\begin{equation*}
\tilde{{\bf h}}=T(g^{-1}){\bf h}(T(g){\bf v},{\bf r},t).
\end{equation*}
Since ${\bf h}$ is dependent on time, then $\tilde{{\bf h}}={\bf h}$. 

We initially applied a shock of minimum amplitudes of $A=0.1, 0.5$, and observed that it was not sufficient to convert 
from a rigidly rotating to a meandering spiral wave. Further, increasing it to $A=1.0$, we could see the conversion 
from rigidly rotating wave to a meandering wave pattern. The conversion was also possible with the shock amplitude 
of $A=1.2$ whereas an increase in the amplitude to $A=1.3$ resulted in the elimination of spiral wave activity, 
thereby causing defibrillation.

\begin{figure*}[t]
\begin{center}
\begin{minipage}[htbp]{0.48\linewidth}
{\centering
\includegraphics[width=0.95\textwidth]{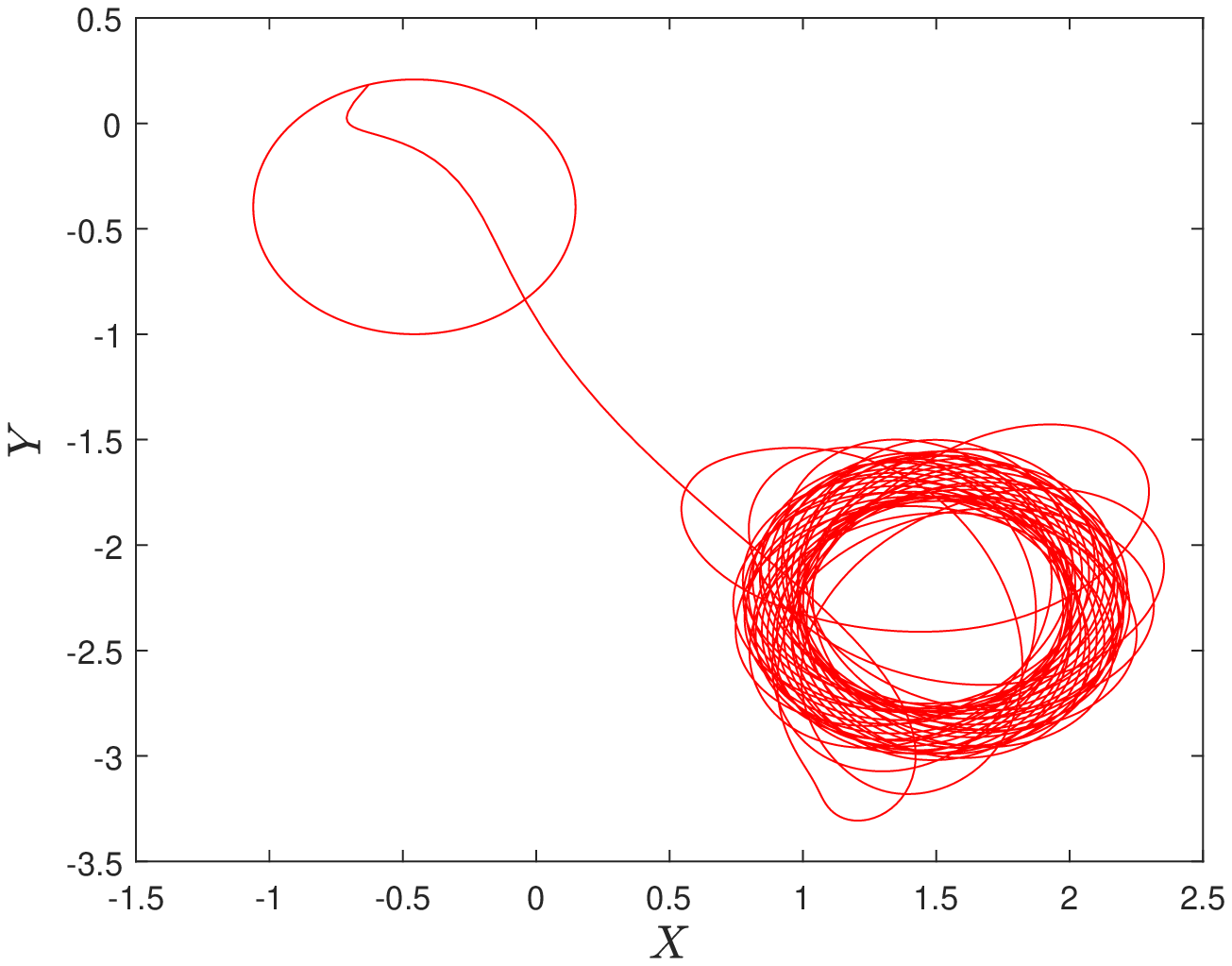}}
\end{minipage}
\begin{minipage}[htbp]{0.48\linewidth}
{\centering
\includegraphics[width=0.95\textwidth]{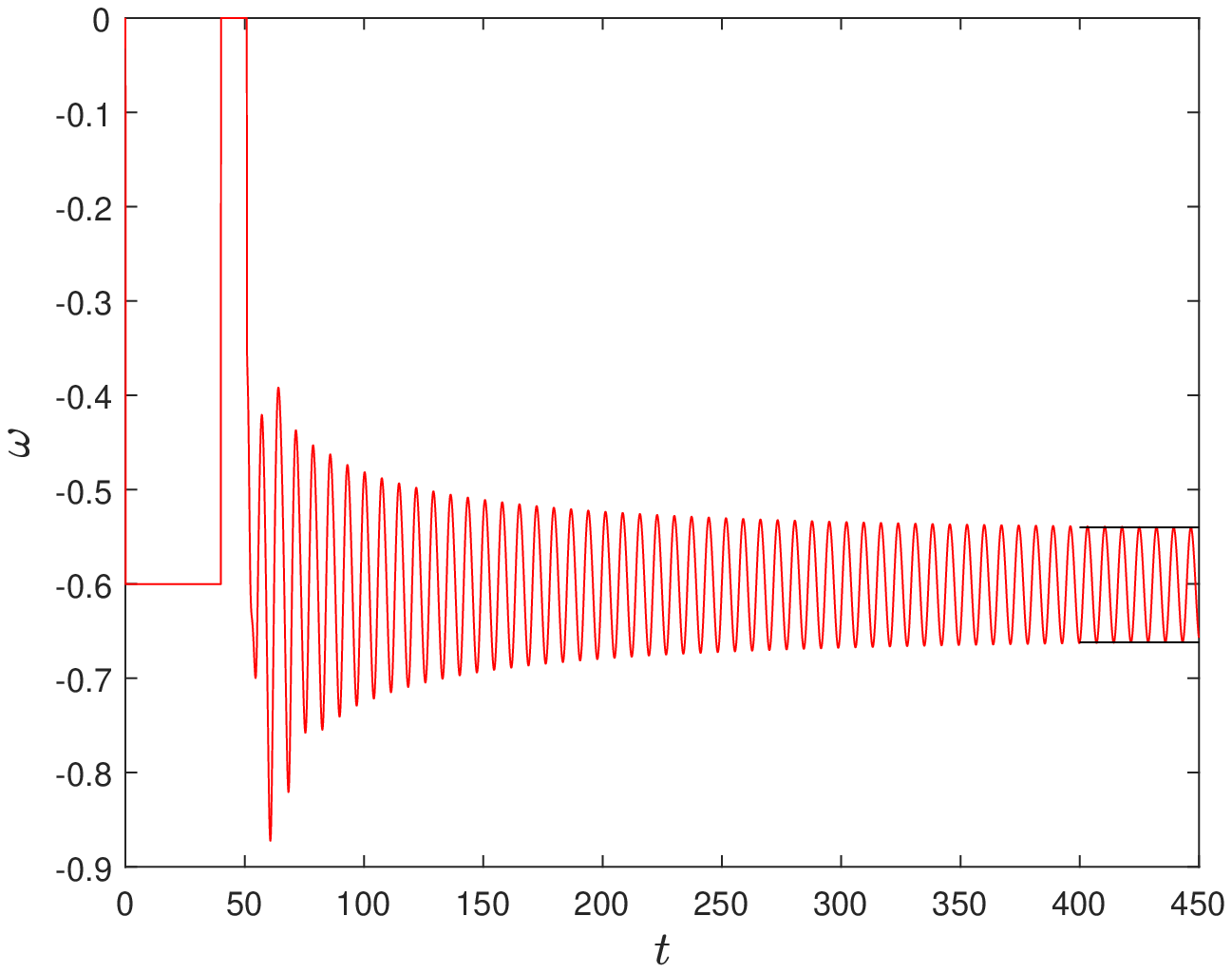}}
\end{minipage}
\begin{minipage}[htbp]{0.48\linewidth}
{\centering
\includegraphics[width=0.95\textwidth]{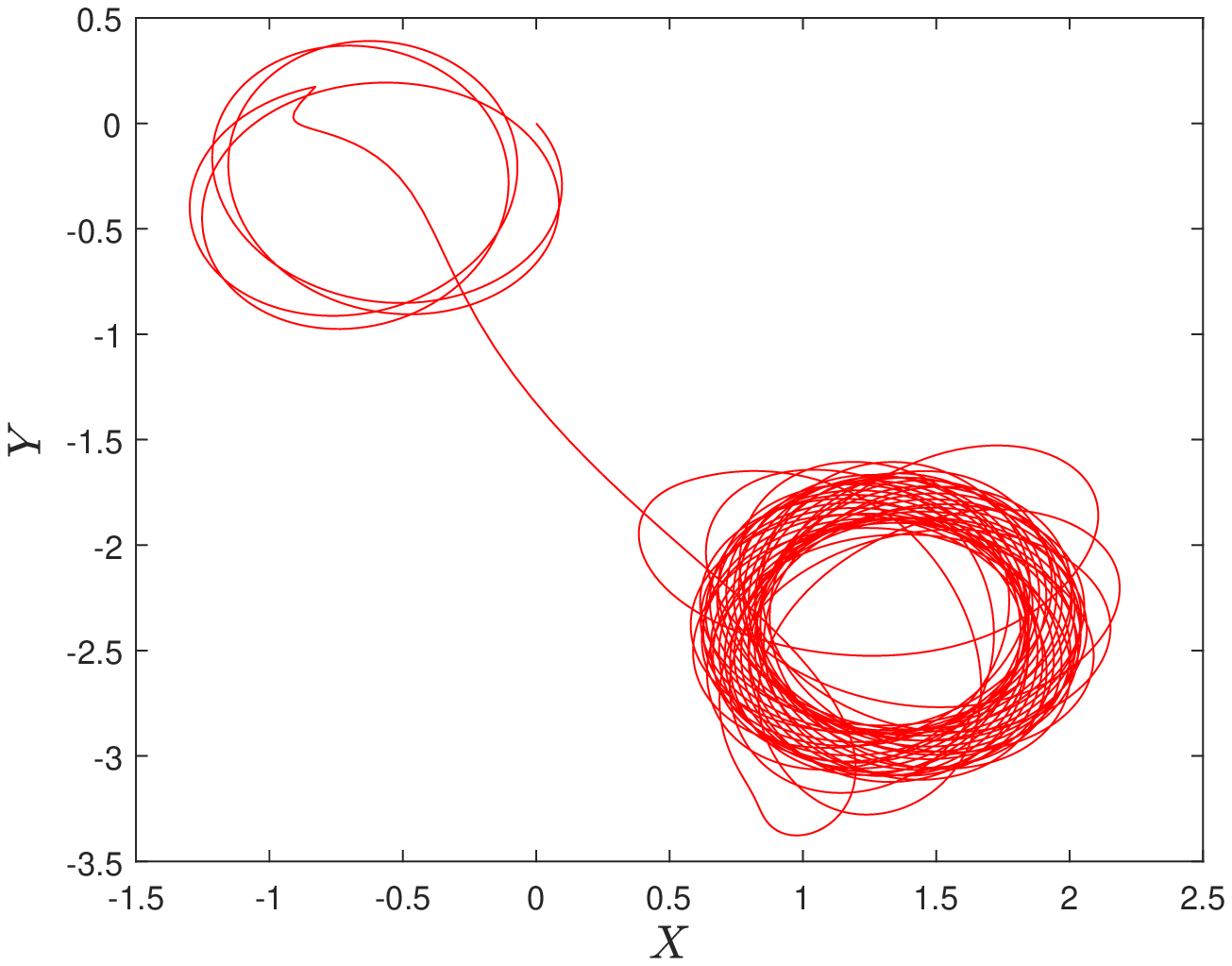}}
\end{minipage}
\begin{minipage}[htbp]{0.48\linewidth}
{\centering
\includegraphics[width=0.95\textwidth]{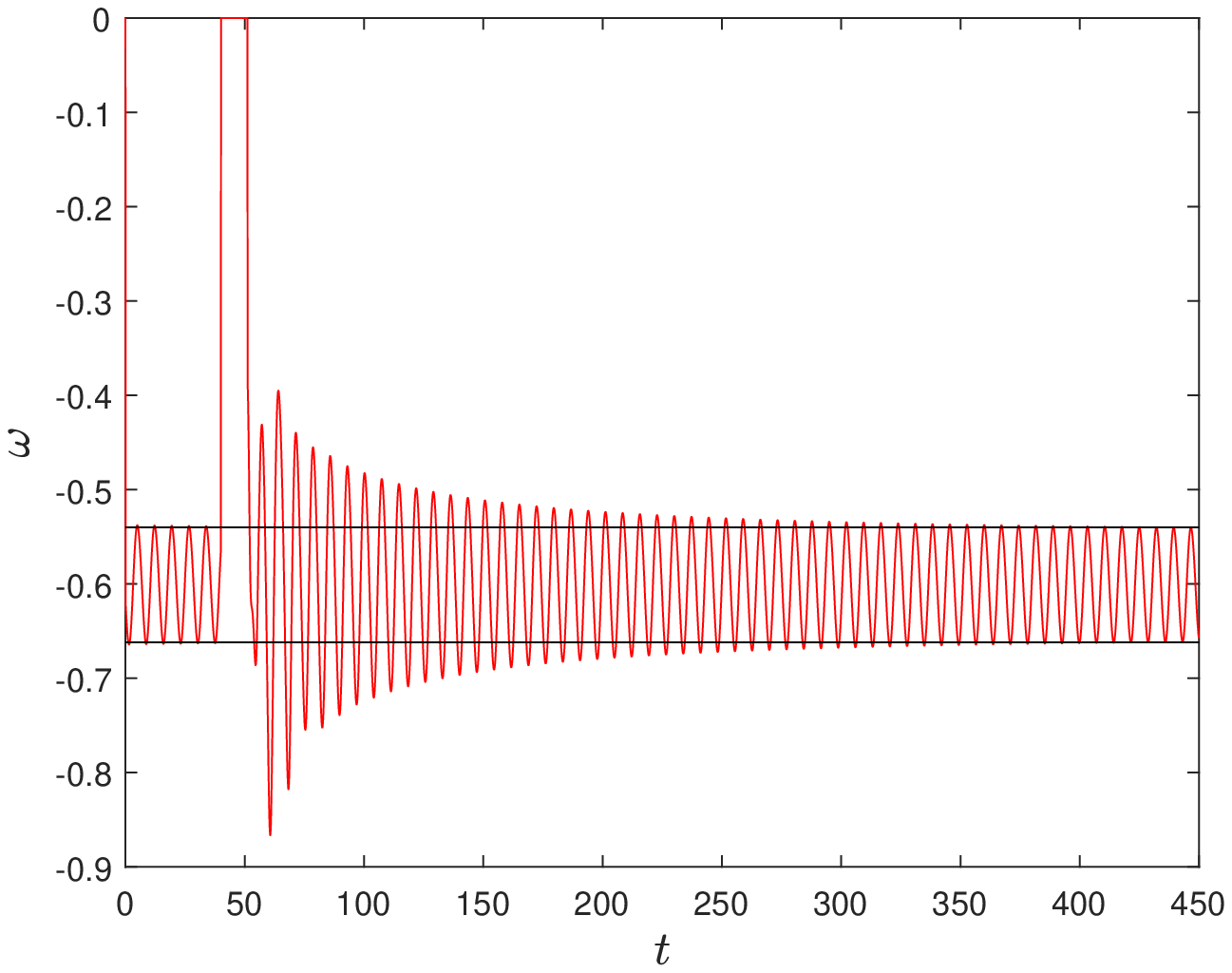}}
\end{minipage}
\caption{(Top) conversion of the rigidly rotating spiral wave to meander; (bottom) attempt to convert a meandering spiral 
wave to another solution. We show both the reconstructed tip trajectory (left) together with the $t$-$\omega$ plots (right),
illustrating the solutions have settled. This is illustrated by using parallel black lines on the $t$-$\omega$ plots.}
\label{fig:shock1}
\end{center}
\end{figure*}

All these shocks were applied after 40,000 steps to check for the conversion. After applying the shock, we can see 
the change in the steady spiral state changes to a periodic state. It concludes that the shock successfully converted 
a rigidly rotating spiral wave solution into a meandering solution, as shown in figure (\ref{fig:shock1}(top)). 

We also tested if we could shock from a periodic solution to a stable solution. We now considered the meandering spiral 
wave solution for $\beta=0.595$ and applied a shock of same amplitude to it. We did not observed the conversion from 
meandering to rigid rotation by a shock of any amplitude. The solution continues to oscillate and converges within 
the same range as before the shock was applied. This implies that conversion is not possible from meander to rigid state. 

Alonso et al.~\cite{negative} noted that the meandering spiral waves with the inward facing petals in two-dimensions 
always have negative filament tension in three-dimensional case, and meandering solutions with outward facing 
petals have positive filament tension. Hence, this gives us an 
explanation for the conversion not being possible in our case, as shown in figure (\ref{fig:shock1}). 
For $\beta=0.595$, we can see that these are the outward meandering spiral waves corresponding to positive tension 
in 3D. In addition, Foulkes et al.~\cite{Foul10} have shown that it is only possible to convert a wave with negative 
filament tension to a wave with positive filament tension. Therefore, for $\beta=0.595$, it is possible to convert 
from rigid rotation to meander but not conversely. 

\section{Convergence analysis}
\begin{figure}[h]
\begin{center}
\begin{minipage}[htbp]{0.48\linewidth}
{\centering
\includegraphics[width=0.95\textwidth]{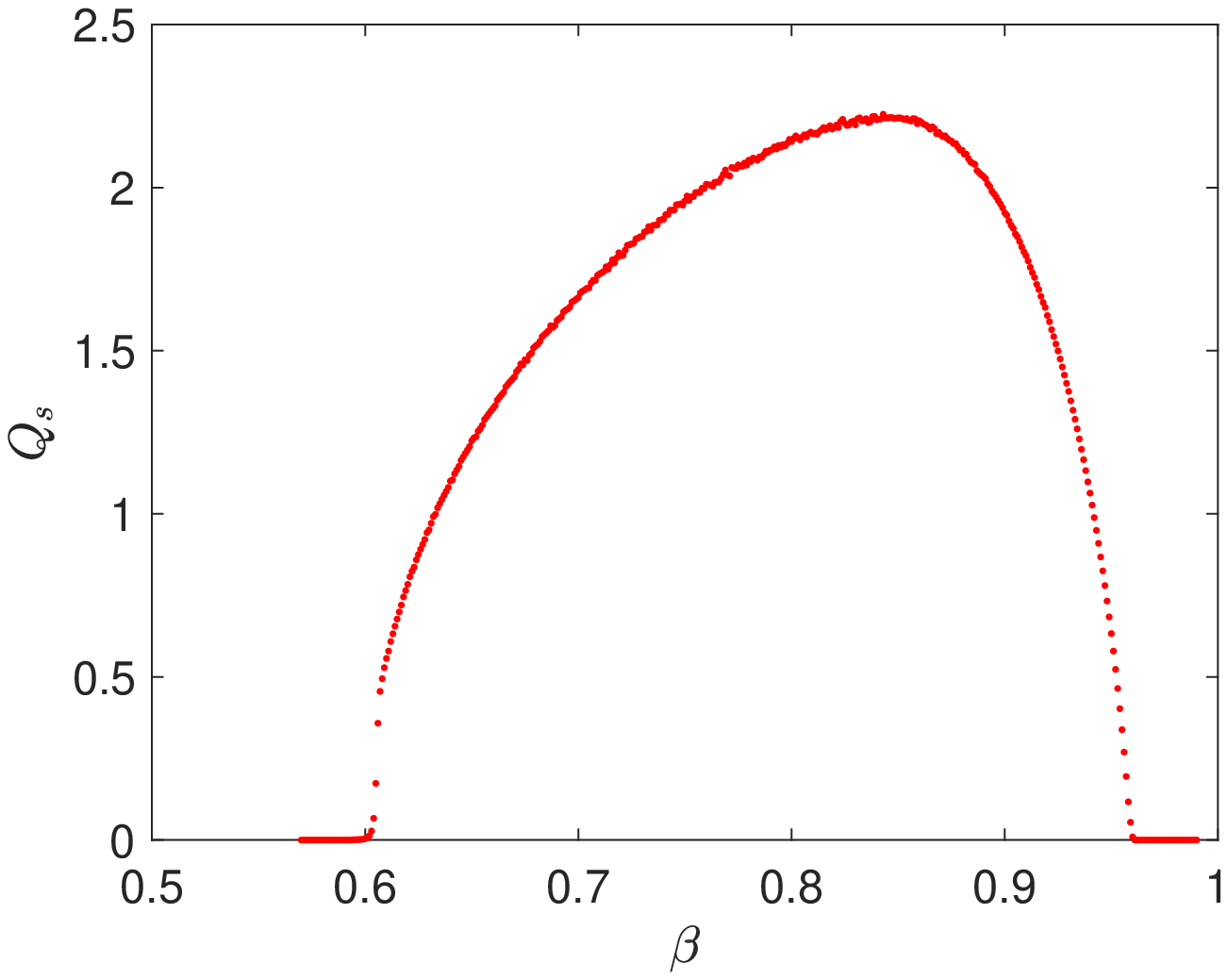}}
\text{($a$)}
\end{minipage}
\begin{minipage}[htbp]{0.48\linewidth}
{\centering
\includegraphics[width=0.95\textwidth]{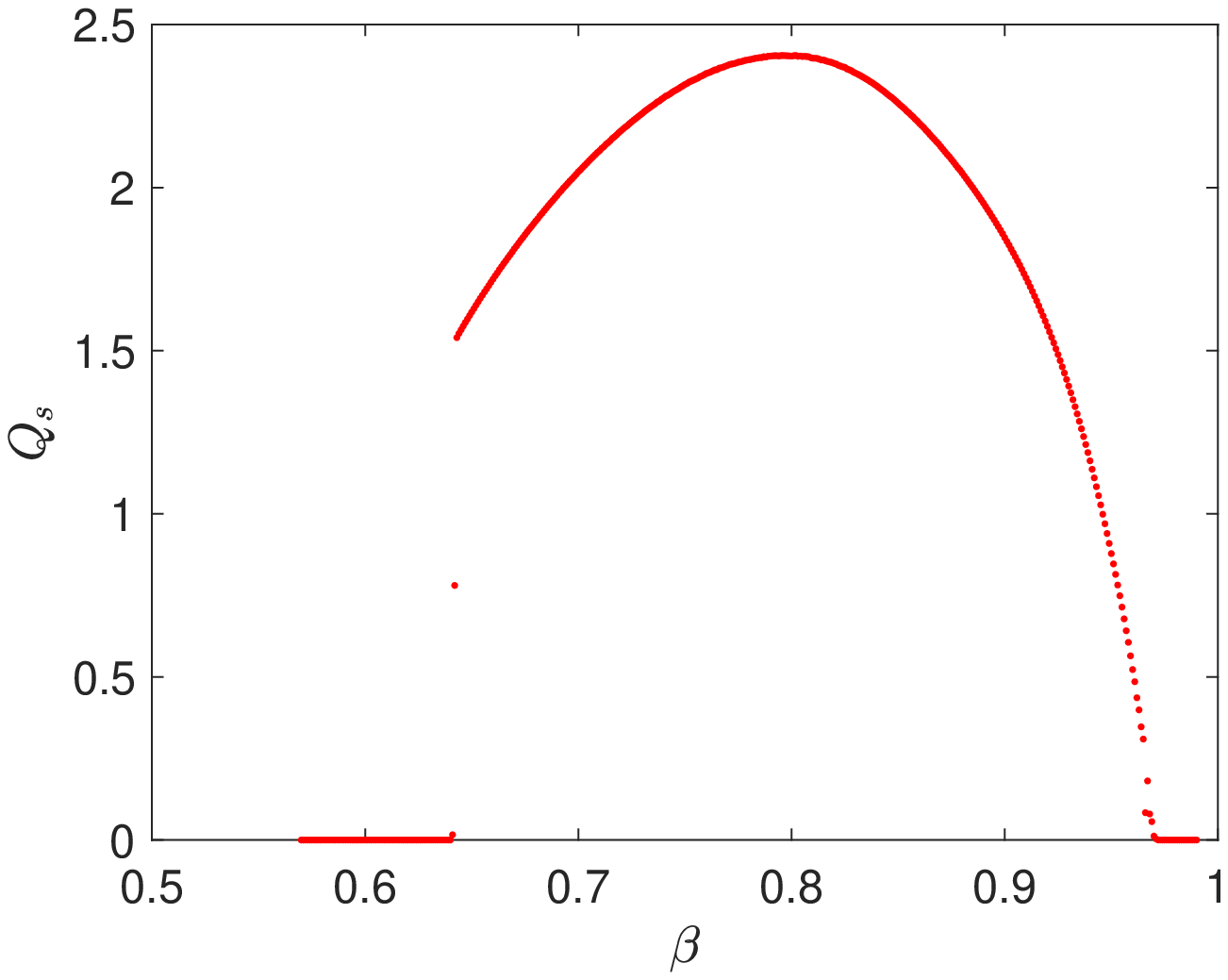}}
\text{($b$)}
\end{minipage}
\begin{minipage}[htbp]{0.48\linewidth}
{\centering
\includegraphics[width=0.95\textwidth]{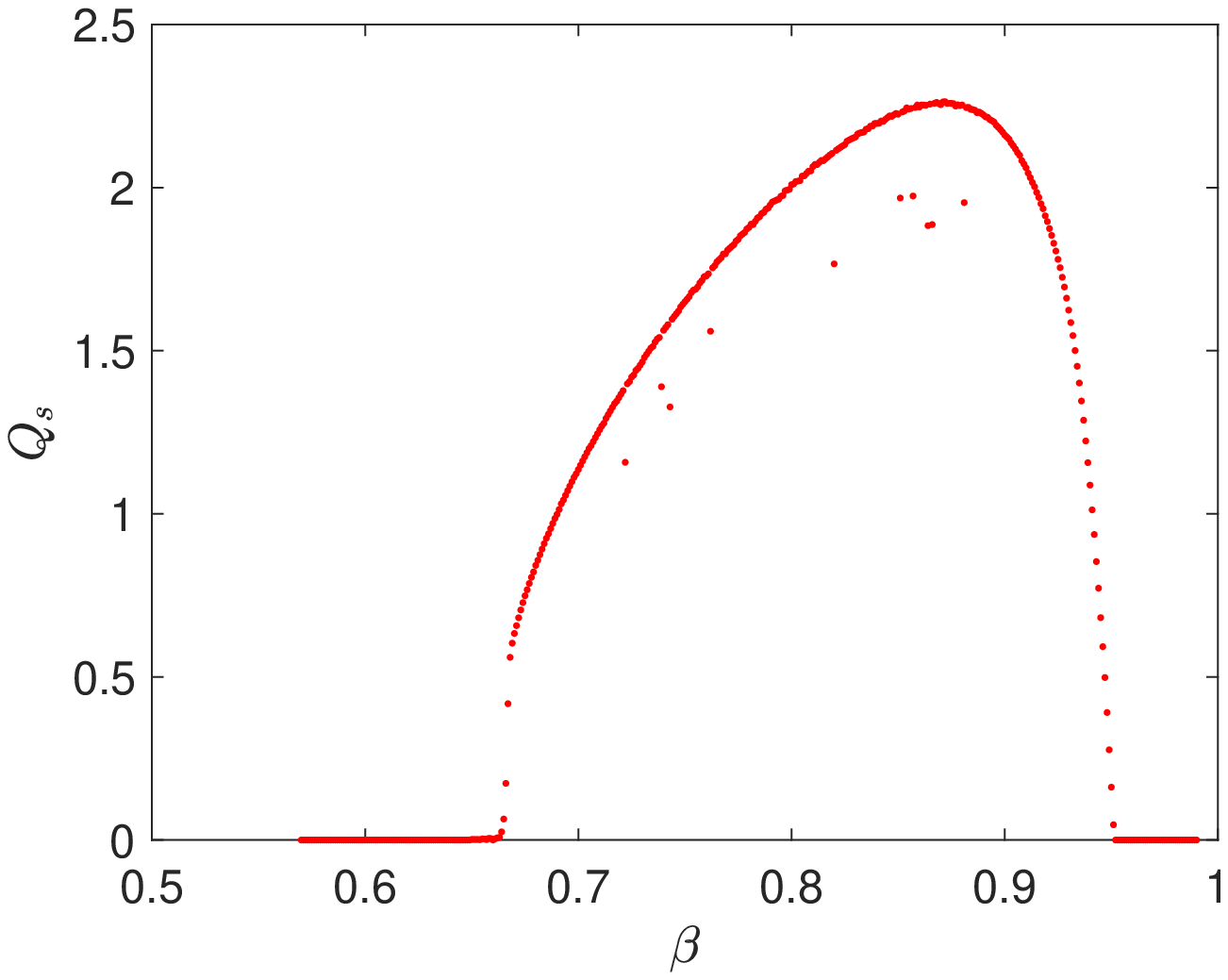}}
\text{($c$)}
\end{minipage}
\begin{minipage}[htbp]{0.48\linewidth}
{\centering
\includegraphics[width=0.95\textwidth]{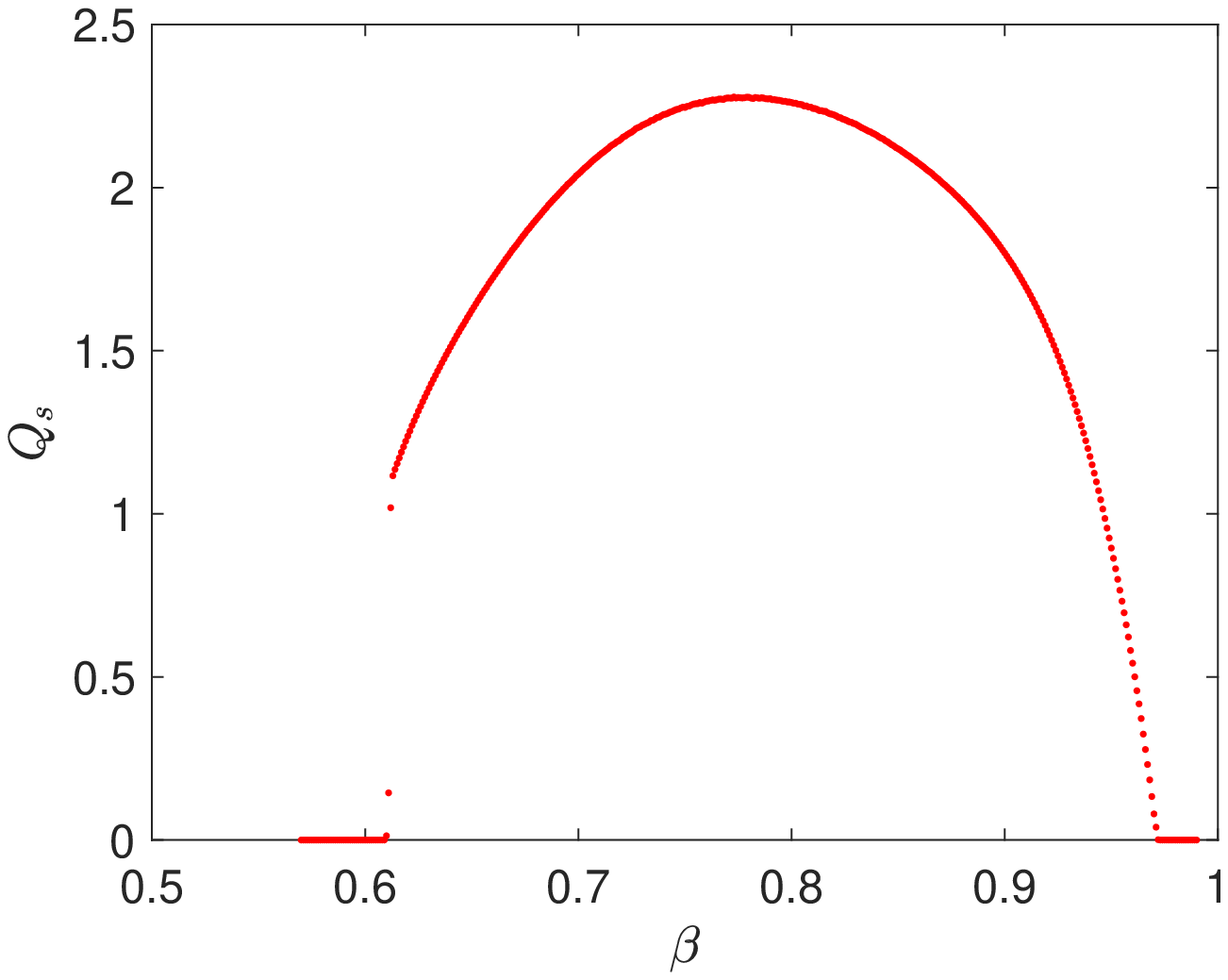}}
\text{($d$)}
\end{minipage}
\begin{minipage}[htbp]{0.48\linewidth}
{\centering
\includegraphics[width=0.95\textwidth]{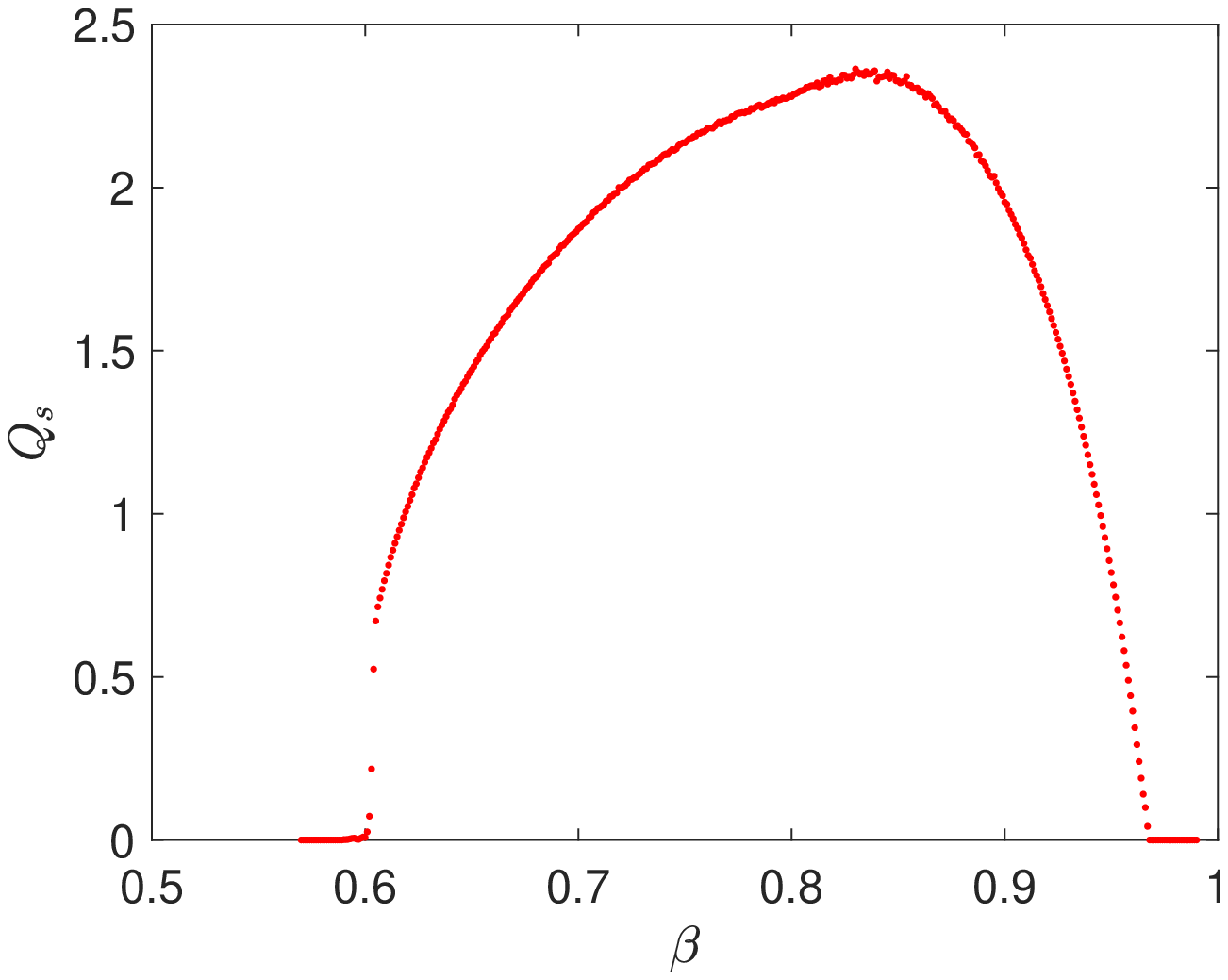}}
\text{($e$)}
\end{minipage}
\begin{minipage}[htbp]{0.48\linewidth}
{\centering
\includegraphics[width=0.95\textwidth]{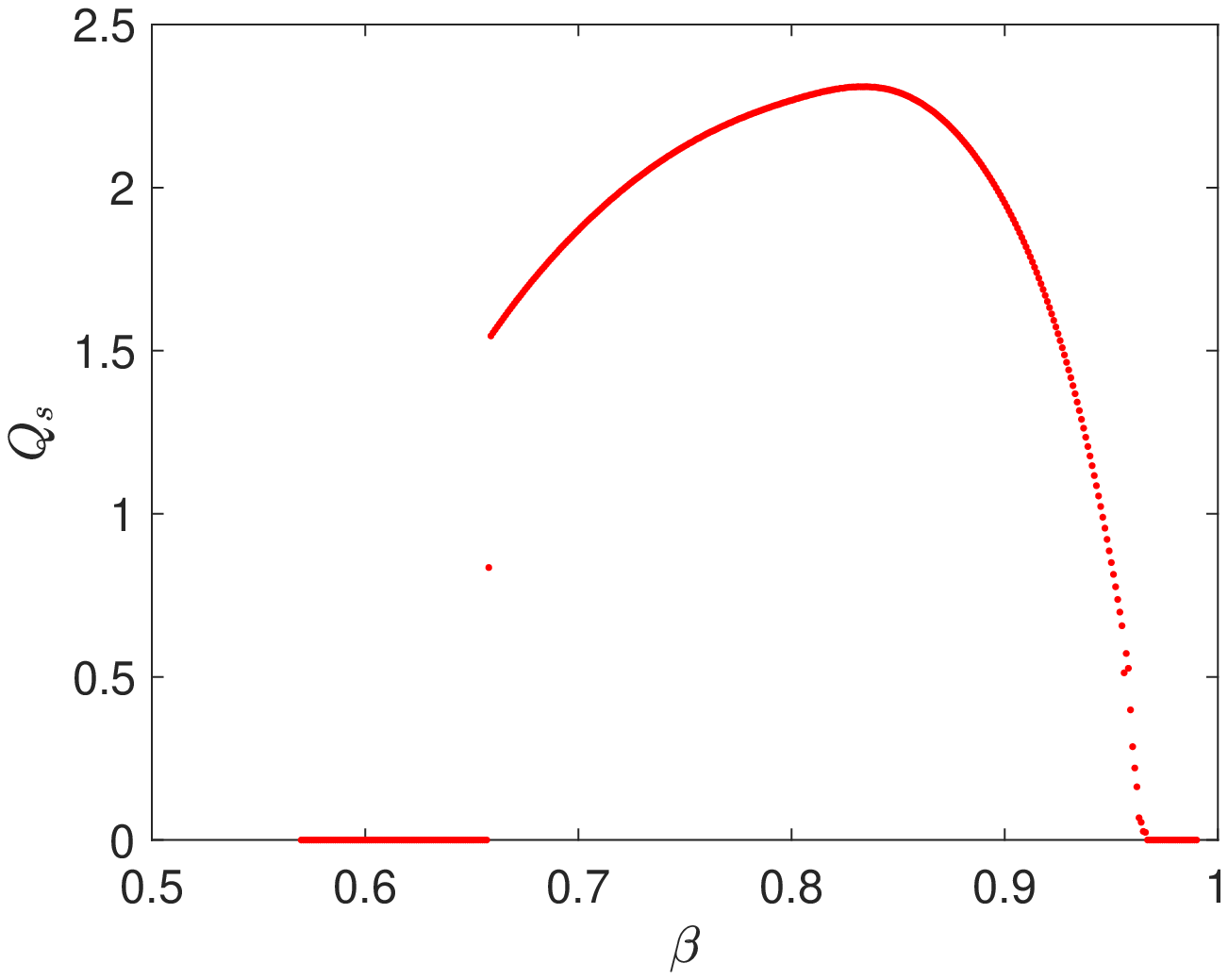}}
\text{($f$)}
\end{minipage}
\begin{minipage}[htbp]{0.48\linewidth}
{\centering
\includegraphics[width=0.95\textwidth]{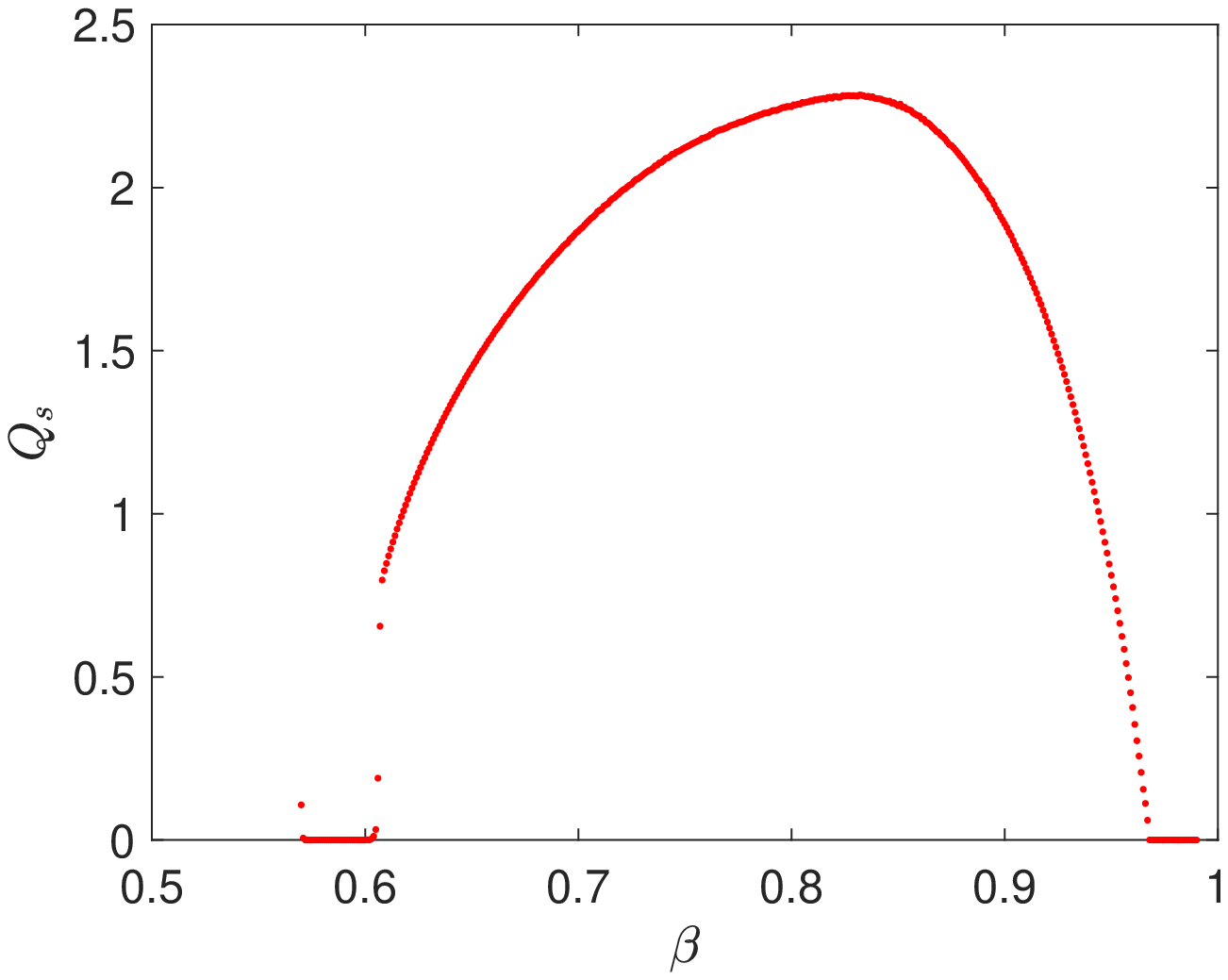}}
\text{($g$)}
\end{minipage}
\begin{minipage}[htbp]{0.48\linewidth}
{\centering
\includegraphics[width=0.95\textwidth]{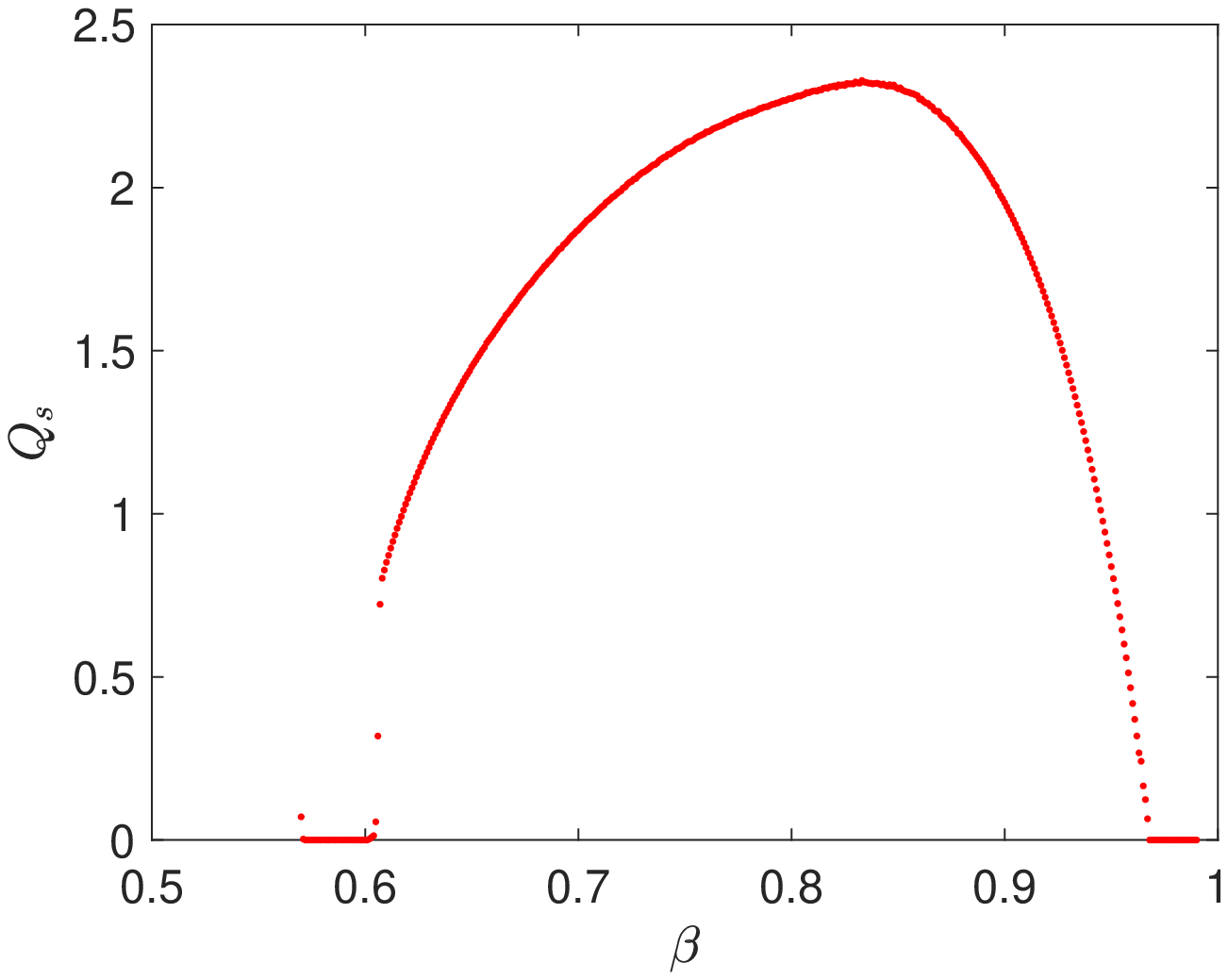}}
\text{($h$)}
\end{minipage}
\caption{Convergence analysis: change in $\Delta_x$ changes $\Delta_t$, keeping $t_s$ constant (a) $\Delta_x=\frac{1}{3}$, (b) $\Delta_x=\frac{1}{8}$; change in $\Delta_x$, allow $t_s$ to change but keeping $\Delta_t$ constant (c) $\Delta_x=\frac{1}{2}$, (d) $\Delta_x=\frac{1}{10}$; change in time step, keeping space step constant (e) $t_s$=0.2, (f) $t_s$=0.025; change in box size (g) $L=15$, (h) $L=60$.} 
\label{fig:conv_fig}
\end{center}
\end{figure}

\seclabel{conv}
We also need to test that the numerical results that we are observed are in fact a true representation of the 
analytical results. We therefore performed several tests to see if the solutions we originally observed converged 
by changing $\Delta_x$, $\Delta_t$ and the domain size, $L$, to prove that it is an accurate representation of the true solution. 
We note that $\Delta_x=L/N_x$ and $\Delta_t=t_s(\Delta_x)^2/4$, so to vary $\Delta_x$ while keeping $L$ fixed, we varied $N_x$.
Similarly, if we vary $\Delta_x$, and keep $t_s$ fixed, then $\Delta_t$ will also vary. If we vary $L$, then in order to 
keep $\Delta_x$ fixed, we must vary $N_x$. A summary of the four tests we conducted is shown below.
\begin{itemize}
\item $\Delta_x$ varied, $t_s$ fixed implying $\Delta_t$ varies, $L$ fixed;
\item $\Delta_x$ varied, $t_s$ varied so that $\Delta_t$ is fixed, $L$ fixed;
\item $\Delta_x$ fixed, $t_s$ varied so that $\Delta_t$ varied, $L$ fixed; and
\item $\Delta_x$, $t_s$ and $\Delta_t$ are fixed, $L$ varies.
\end{itemize}
The default values of the various numerical parameters are $L=30$, $\Delta_x = 0.2$, $\Delta_t = 0.001$, $t_s = 0.1$. 
Figure (\ref{fig:conv_fig}) shows a selection of the convergence tests that we conducted, plotting $\beta$ 
against $Q_s$. 

The first test, shown in figure (\ref{fig:conv_fig}(a)-(b)) was to vary $\Delta_x$ while keeping $t_s$ 
fixed, meaning that $\Delta_t$ varies too. It was observed that
with the decrease in the space step, the size of the left hand hysteresis region becomes more prominent. 
We also note that other changes included an increase in the $Q_s$ as a function of $\beta$, and the bifurcation
points also changed. However, despite all these changes, the overall shape of the graph depicting the growth of 
quotient size against $\beta$ remains the same.

The next test varied the $\Delta_x$ but this time ensured that $\Delta_t$ was fixed by varying $t_s$. Partial 
results are shown in figure (\ref{fig:conv_fig}(c)-(d)). As with the first test, certain features such as the 
location of the bifurcation points, the maximum value of $Q_s$, etc, changed as we varied $\Delta_x$, but the 
overall shape of the arc remained qualitatively the same.

The third test varied $\Delta_t$ while keeping all other numerical parameters fixed, figure (\ref{fig:conv_fig}(e)-(f)). 
While we ran the simulations for different time steps, instabilities were present for the larger time steps but these 
instabilities tend to decrease in size and frequency as we decreased the time step. But, yet again, while features such
as the location of the bifurcation points, maximum value of $Q_s$, etc changed slightly, the overall 
shape of the curve stayed the same.

The final test is the convergence in box size and some of the results are shown in figure (\ref{fig:conv_fig}(g)-(h)). 
We considered two variants: first with half the box size and the other one when we double the original box size.
In both the cases, we got almost the same Hopf bifurcation points. Only the properties (area, $\beta$ with 
maximum quotient size) affecting quotient size differs while we increase the box size. 

These tests, show that as we vary the numerical parameters, we change some of the features of the results such as 
the location of the Hopf bifurcation points, the maximum value of $Q_s$, but one thing remained clear, and that
was the overall shape of the curve remained similar throughout the studies. It is clearly evident that the 
original data represents a true likeness of the actual data.  

\section{Discussion}\seclabel{disc}
We have examined the type of bifurcation responsible for the transition from rigidly rotating spiral wave solutions to meandering spiral wave solutions using FHN model. 

We chose FHN model, where $\beta$ was considered as a bifurcation parameter. It is also important to note that we 
conducted numerical simulations in the FOR moving with the tip of a spiral wave. This helped us to overcome the 
drawbacks of previous studies which were unable to quickly study the large core spirals~\cite{Bark94a}. We were 
able to successfully study the bifurcation analysis of spiral waves for the full range of $\beta$ depicting the 
transition from rigid to meander and then back to rigid motion. 

Since, our analysis was conducted in comoving FOR, we study the limit cycle solutions within meandering spiral 
waves described by the advection coefficients. Here, the limit cycles are not of any particular shape. Therefore, 
we presented a new technique to check for the type of bifurcation responsible for the transition from rigid to 
meander motion in FHN model with a specific set of parameters. We calculated the quotient size of limit cycles 
and plotted it as a function of $\beta$. With this approach, we discovered that unlike other 
authors~\cite{Bark94a, Swinney}, there is a subcritical Hopf bifurcation responsible for the transition from 
rigid rotation (RW) to meander states (MRW) in the parameter space that we studied. Near the Hopf bifurcation points, it was observed 
that there is a discontinuous jump, being more prominent on the left hand side of the bifurcation diagram. Due to 
this discontinuous growth in $Q_s$, the analysis was conducted in a reverse direction across the $\beta$-range, 
which disclosed a region of hysteresis where one could find two spiral wave solutions for a single set of parameters. 
It made it more clear about the existence of subcritical bifurcations within FHN model. 

The hysteresis region is the region of bistability where there exists two types of solutions for all the values 
of $\beta$ within the region. We investigated this region by applying the technique of single shock defibrillation. 
It us a result that only a rigidly rotating spiral wave can be converted into a meandering spiral wave but not conversely.

Furthermore, the $\beta$-step in the system was arbitrarily closely approximated by other $\beta$-steps with 
distinguishably different bifurcation points, resulting in wider hysteresis region. In other words, we can say 
that a small change to initial conditions may lead to different behaviour, thereby showing the sensitive 
dependence on initial conditions. 

We have also demonstrated the numerical convergence analysis which confirms the true representation of the 
solutions. It highlighted that the overall shape of the bifurcation diagrams remains the same for all numerical 
and model parameters. It was clearly noted that the change in the numerical or model parameters does not 
qualitatively affect the bifurcation diagram depicting the growth of quotient size. Adding to it, these 
tests can be used in future to check for different numerical parameters which would enable us to run faster 
simulations for various studies.

\begin{acknowledgments}
This study has been funded in part by a Vice-Chancellor's PhD Scholarship
from Liverpool Hope University. The authors are grateful for discussions 
with Prof. Vadim Biktashev (University of Exeter, UK) and Dr. Irina
Biktasheva (University of Liverpool, UK). 
\end{acknowledgments}

\end{document}